\documentclass[moor,nonblindrev]{informs1}

\usepackage{physics}
\usepackage{mathtools}
\usepackage{tensor}
\usepackage{comment}

\usepackage[mathscr]{euscript}   

\newcommand{\half}{\tfrac{1}{2}}
\newcommand{\asymx}{\mathop{\sim}}
\newcommand{\asym}[1]{\mathrel{\asymx_{#1}}}

\newcommand{\Rez}{\mathop{\rm Re} \nolimits}

\newcommand{\sinc}{\mathop{\rm sinc} \nolimits}

\newmuskip\pFqmuskip
\newcommand*\pFq[6][8]{%
  \begingroup 
  \pFqmuskip=#1mu\relax
  \mathcode`\,=\string"8000
  \begingroup\lccode`\~=`\
  \lowercase{\endgroup\let~}\pFqcomma
  {}_{#2}F_{#3}{\left[\genfrac..{0pt}{}{#4}{#5};#6\right]}%
  \endgroup
}
\newcommand{\pFqcomma}{\mskip\pFqmuskip}

\DeclareMathOperator{\acosh}{acosh}

\def\wscl{0.8}
\def\hscl{0.56}

\usepackage{natbib}
\NatBibNumeric
\def\bibfont{\small}%
\bibpunct[, ]{[}{]}{,}{n}{}{,}%

\definecolor{myblue}{RGB}{0, 20, 114}
\usepackage[colorlinks=true,breaklinks=true,bookmarks=true,urlcolor=myblue,
     citecolor=myblue,linkcolor=myblue,bookmarksopen=false,draft=false]{hyperref}

\def\EMAIL#1{\href{mailto:#1}{#1}}

\usepackage{algorithmic}
\usepackage{algorithm}

\TheoremsNumberedThrough     

\EquationsNumberedThrough    

\MANUSCRIPTNO{0} 

\def\theARTICLETOP{}
\def\theLRHSecondLine{}
\def\theRRHSecondLine{}

\begin{document}


\RUNAUTHOR{Zuk and Kirszenblat}

\RUNTITLE{Non-Preemptive Two-Level Priority Queue}

\TITLE{Analytic Approach to the Non-Preemptive Markovian Priority Queue}

\ARTICLEAUTHORS{%
\AUTHOR{Josef Zuk}
\AFF{Defence Science and Technology Group, Melbourne, Australia,
\EMAIL{josef.zuk@defence.gov.au}}
\AUTHOR{David Kirszenblat}
\AFF{Defence Science and Technology Group, Melbourne, Australia,
\EMAIL{david.kirszenblat@defence.gov.au}}
} 

\ABSTRACT{%
A new approach is developed for the joint queue-length distribution
of the two-level non-preemptive Markovian priority queue
that allows explicit and exact results to be obtained.
Marginal distributions are derived for the general multi-level problem.
The results are based on a representation of the joint queue-length
probability mass function as a single-variable complex contour integral,
that reduces to a real integral on a finite interval arising from a cut on the real axis.
Both numerical quadrature rules and exact finite sums, involving Legendre polynomials
and their generalization,
are presented for the joint and marginal distributions.
A high level of accuracy is demonstrated across the entire ergodic region.
Relationships are established with the waiting-time distributions.
Asymptotic behaviour in the large queue-length regime is extracted.
}%

\MSCCLASS{Primary: 90B22; secondary: 60K25, 60J74}
\ORMSCLASS{Primary: Queues: Priority; secondary: Queues: Markovian}

\HISTORY{Date created: September 05, 2023. Last update: December 07, 2023.}

\maketitle

%

\section{Introduction}
\label{intro}
This work is concerned with both algebraically exact and highly accurate numerical evaluations
of the joint and marginal queue length distributions for the non-preemptive
Markovian priority queue, based on a representation of the joint queue-length
probability mass function (PMF) as a single-variable complex contour integral,
making use of its pole and cut singularities on the real axis.

There exist numerous real-world applications of the non-preemptive priority queue,
including telecommunications \citep{NP:Cohen56},
health care \citep{NP:Almehdawe13,NP:Hou20,NP:Taylor80},
beam scheduling in phased-array radars \citep{NP:Orman95,NP:Orman96}, and
air traffic control \citep{NP:Pestalozzi64}, to name just a few.
Here, we revisit the two-level problem from a different perspective
to that considered by us in \citep{NP:Zuk23} and the previous literature
cited therein.
The two-level case is amenable to a number of exact analytical approaches
that are not easily accessible to the general multi-level problem.
Furthermore, the two-level problem occupies a distinguished position,
because knowledge of the low-priority marginal queue-length distribution for
the two-level case is sufficient for the evaluation of all the marginal queue-length
distributions in the general problem.
For this reason, the two-level problem justifies some special attention.

The two-level problem was recently discussed in \citep{NP:Zuk23},
where a detailed survey of previous work was also given.
The most recent work on the problem, prior to our own, appeared in \citep{NP:Sleptchenko03}.
The multi-level problem for an arbitrary number of priority levels was
solved in \citep{NP:Zuk23A}.
In that work, an explicit exact expression was derived for the multi-variate
joint queue-length probability generating function (PGF),
but the extraction of the joint and marginal queue-length
PMFs from the PGF was numerical, rather than analytical.
Here, we describe an analytical approach to the two-level problem.
It is, at present, not clear whether an analogous analytic method could be easily
extended to the general multi-level problem.

In the present approach, both the marginal and joint PMFs are represented as
one-dimensional real-valued integrals over a finite interval.
The integral representation is amenable to efficient numerical evaluation
via Gaussian quadrature in an iterative form that allows a desired level of accuracy
to be set {\it a priori}.
It is equally amenable to an exact evaluation as a finite sum of easily
computable functions that generalize the familiar associated Legendre functions
in a novel way.
Both methods produce highly accurate results for parameter values spanning the entire
ergodic region, including values of the total traffic intensity close to unity.

The analytical structure of the integrals arising here is very similar to that
encountered for the partial busy period of an M/M/$c$ queue \cite{NP:Karlin58,NP:Zuk23B},
and for the waiting-time marginal distributions in \citep{NP:Davis66}.
In fact, we are able to provide alternative derivations of the
waiting-time marginal distributions,
either from the queue-length marginal PGFs, or from the queue-length marginal PMFs.
In the latter case, we find a new representation as a convergent Laguerre series,
where the Laguerre polynomials carry the time dependence, while the sequence
of coefficients is a binomial transform of the queue-length PMF.

For multi-class non-priority queueing models,
\citet{NP:Bertsimas97} derive a representation of the multi-variate joint PGF
for the queue-length distribution as a multi-dimensional time integral over all the
waiting-time marginals, based on the distributional form of Little's law.
However, this approach does not extend to the priority queue because the
underlying no-overtaking assumption is violated.
Consequently, in the present problem, the joint queue-length distribution cannot
be derived from the waiting-time marginals.
Here, and in \citep{NP:Zuk23A}, a direct derivation of the joint queue-length
distribution is provided.

The integral representations of the marginal and joint PMFs facilitate the
extraction of the asymptotic behaviour in the large queue-length limits;
something that appears to be difficult to achieve from other starting points.
The integral representations are also robust to the inclusion of various
modelling complications.
An example arises in the application to the ambulance ramping problem \citep{NP:Almehdawe13}.
Here, the priority levels are determined by patient acuity levels,
but patients also fall into separate arrival classes depending on whether they
are entering the hospital emergency department directly from an ambulance, from
an ambulance offload zone or as walk-ins. Each priority level may span a number of
arrival classes, and each arrival class has a distinct mean arrival rate.
This results in a multi-level, multi-class queueing model
with priorities and classes entwined.
The queue-length PMFs acquire a somewhat more complex integral representation,
but the simple quadrature rules derived for the present model remain applicable
and efficient. We shall make use of this in a forthcoming paper.

The remainder of the paper is organized as follows: In Section~\ref{PGF},
we introduce an approach to the problem based on a PGF
of a single complex variable that sums the joint distribution only over the low-priority
indices ({\it i.e.}\ queue lengths), and we elucidate the analytic structure of the PGF
in the complex plane.
In Section~\ref{Marginal}, we derive an integral representation for the low-priority
marginal queue-length distribution, noting that the high-priority marginal is trivially geometric,
and proceed to establish iterative quadrature rules.
We show that the results are related to one of Appell's bivariate hypergeometric functions.
The asymptotic behaviour in the limit
of large queue length is also derived directly from the integral representation.
In Section~\ref{WaitTime}, we discuss the marginal waiting-time distributions
for the general multi-level non-preemptive priority problem.
As a by-product of our results on the queue-length marginal distributions, we are able to invoke
the distributional form of Little's law to reproduce a previously established integral
form for the waiting-time marginals, which may be related to a bivariate confluent
hypergeometric function of Humbert,
as well as deriving a new representation as a Laguerre series.
In Section~\ref{Joint}, we derive an integral representation for the joint
queue-length distribution, and establish the relevant iterative quadrature rules.
In Section~\ref{Exact}, we first show that the low-priority marginal PMF can be
expressed as a finite sum of Legendre polynomials.
For the joint PMF, an exact representation as a finite sum from a family of functions $Q_m^\ell(\chi)$
is derived. These functions reduce to the Legendre polynomials in the special case
\mbox{$Q_m^0(\chi) = P_m(\chi)$},
and are closely related to the associated Legendre functions when
\mbox{$\abs{\ell} \leq m$},
but are defined for all integral values of $\ell$.
Various properties of these functions are elucidated.
The section concludes with an efficient and robust numerical implementation of the exact expressions.
In Section~\ref{Numerical}, we introduce a number of measures of performance (MOP) to study the
numerical behaviour of the two approaches developed in the present paper, testing them
against one other, against known exact relationships, and comparing with two other algorithms
previously derived for this problem \cite{NP:Zuk23}.
Some concluding remarks follow, and
various technical details appear in the appendices.

\section{Probability Generating Function}
\label{PGF}
As in \citep{NP:Zuk23}, we consider a non-preemptive priority queue with
$K$ priority levels, each with a distinct
Poisson arrival rate $\lambda_k$,
\mbox{$k = 1,2,\ldots,K$}
and corresponding level traffic intensity\footnote{Consistent with \citep{NP:Gail88,NP:Kao90},
  the notation $\rho_k$ reserved for $\rho_k \equiv \lambda_k/\mu$, so that $r_k = \rho_k/c$.}
\mbox{$r_k = \lambda_k/(c\mu)$},
leading to a total traffic intensity for the aggregation of all arrivals of
\mbox{$r = \sum_{k=1}^K r_k$}.
We adopt the usual convention that smaller priority-level indices $k$
represent higher priorities. Thus, $r_1$ denotes the traffic intensity
associated with the highest priority level.
We assume a common exponential service rate $\mu$
among all priority levels.
The number of servers is arbitrary and is denoted by $c$.
For the discussion of the joint queue-distribution, we limit the number of priority levels to
\mbox{$K = 2$}, and adopt the notation
\mbox{$r_{\text{hi}} \equiv r_1$},
\mbox{$r_{\text{lo}} \equiv r_2$},
so that
\mbox{$r= r_{\text{hi}} + r_{\text{lo}}$}.
We also denote the fraction of high-priority arrivals (abbreviated `hifrac') by
\mbox{$\nu \equiv r_{\text{hi}}/r$},
so that
\mbox{$0 \leq \nu \leq 1$}.
This allows us to study all possibilities for a given total traffic intensity $r$ by
independently varying $\nu$ over its full range.

Let
\mbox{$P_{\text{full}}(\ell,m)$}
denote the unconditional queue-length PMF with $\ell$ high-priority
clients and $m$ low-priority clients in the queue.
Let
\mbox{$P(\ell,m)$}
denote the same queue-length PMF conditioned on all servers being busy,
to which we shall refer as the wait-conditional PMF.
These two quantities are related by
\begin{equation}
P_{\text{full}}(\ell,m) = P_{\text{NW}}\delta_{\ell 0}\delta_{m0}
     + (1 - P_{\text{NW}})P(\ell,m) \;,
\end{equation}
where $P_{\text{NW}}$ is the no-wait (or delay) probability.
An expression for its value is derived in Appendix~\ref{NoWait}.
The wait-conditional PMF does not directly depend on the
number of servers, but only indirectly through its dependence on the
total traffic intensity.
Therefore, the model for
\mbox{$P(\ell,m)$}
depends on just two parameters, namely, $r$ and $\nu$,
that may be varied independently in the intervals
\mbox{$0 \leq r < 1$},
\mbox{$0 \leq \nu \leq 1$}.
\citet{NP:Cohen56} introduced a PGF with a discrete index
for the unconditional PMF given by
\begin{equation}
G_\ell(z) \equiv \sum_{m=0}^\infty z^m P_{\text{full}}(\ell,m) \;.
\end{equation}
The corresponding wait-conditional PGF
\mbox{$g_\ell(z)$}
is related to the full PGF by
\begin{equation}
G_\ell(z) =  P_{\text{NW}}\delta_{\ell 0}
     + (1 - P_{\text{NW}})g_\ell(z) \;.
\end{equation}

To obtain an expression for
\mbox{$g_\ell(z)$},
we introduce
\mbox{$\zeta \equiv \zeta_\pm(z)$}
as the two branches that solve the $z$-dependent quadratic equation
\begin{equation}
\zeta^2 - (1 + r - r_{\text{lo}}z)\zeta + r_{\text{hi}} = 0 \;.
\end{equation}
Thus, we can write
\begin{equation}
\zeta_\pm(z) = \half\left[b(z) \pm \sqrt{\Delta(z)}\right] \;,
\end{equation}
where
\begin{equation}
\Delta(z) \equiv b^2(z) - 4r_{\text{hi}} \;, \quad
     b(z) \equiv 1 + r - r_{\text{lo}}z \;,
\end{equation}
and it is useful to note the quadratic identities
\begin{equation}
\zeta_+(z) + \zeta_-(z) = b(z) \;, \quad
     \zeta_+(z){\cdot}\zeta_-(z) = r_{\text{hi}} \;,
\label{SumProd}
\end{equation}
which also imply that
\begin{equation}
[1 - \zeta_+(z)]{\cdot}[1 - \zeta_-(z)] = r_{\text{lo}}(z-1) \;.
\end{equation}
Notable special values are
\mbox{$\zeta_+(1) = 1$},
\mbox{$\zeta_-(1) = r_{\text{hi}}$}.
Following \citet{NP:Cohen56}, the wait-conditional PGF is given by
\begin{equation}
g_\ell(z) = g_{\text{lo}}(z){\cdot}\left[1 - \zeta_-(z)\right]\zeta_-^\ell(z) \;, \quad
     g_{\text{lo}}(z) = \frac{1-r}{\zeta_+(z) - r} \;.
\label{glz}
\end{equation}
The quantity
\mbox{$g_{\text{lo}}(z)$}
represents the PGF for the low-priority marginal distribution, as can be seen from
\begin{equation}
\sum_{\ell=0}^\infty g_\ell(z) = g_{\text{lo}}(z) \;.
\end{equation}
These results also follow as a special case of the multi-variate joint PGF
for the general multi-level non-preemptive queue derived in \citep{NP:Zuk23B}.

The discriminant
\mbox{$\Delta(z)$}
can be factorized according to
\begin{equation}
\Delta(z) = r_{\text{lo}}^2(z - x_+)(z - x_-) \;,
\end{equation}
with real strictly positive roots
\begin{equation}
x_\pm = 1 + (1 \pm \sqrt{r_{\text{hi}}})^2/r_{\text{lo}} \;.
\end{equation}
It is useful to introduce their difference
\begin{equation}
    x_{\text{dif}} \equiv x_+ - x_- = 4\sqrt{r_{\text{hi}}}/r_{\text{lo}} \;.
\end{equation}
One should observe that
\mbox{$\Delta(x) < 0$}
for
\mbox{$x_- < x < x_+$},
which implies the presence of a cut on this finite real interval,
and  that
\mbox{$\Delta(1/r) = (r^2 - r_{\text{hi}})^2/r^2 > 0$}.
Thus, the pole at $x = 1/r$
does not lie in the open cut interval
\mbox{$(x_-,x_+)$}.

\section{Marginal Distributions}
\label{Marginal}
The high and low priority marginal PMFs are given, respectively, by
\begin{equation}
P_{\text{hi}}(\ell) = \sum_{m = 0}^\infty P(\ell,m) \;, \quad
    P_{\text{lo}}(m) = \sum_{\ell = 0}^\infty P(\ell,m) \;.
\end{equation}
An immediate consequence of (\ref{glz}) is the simple geometric form
for the high-priority marginal
\begin{equation}
P_{\text{hi}}(\ell) = g_\ell(1) = (1-r_{\text{hi}})r_{\text{hi}}^\ell \;.
\end{equation}
Thus, only the low-priority marginal
\mbox{$p_n \equiv P_{\text{lo}}(n)$}
is of interest for further analysis.
By Cauchy's integral theorem, we have
\begin{equation}
p_n = \oint_\mathcal{C}\frac{dz}{2\pi i}\, \frac{1}{z^{n+1}} g_{\text{lo}}(z) \;,
\end{equation}
where the integration contour $\mathcal{C}$ is taken to be an anti-clockwise circle
centred on the origin, of radius $\eta$ within the radius of convergence of the function
\mbox{$g_{\text{lo}}(z)$},
namely
\mbox{$0 < \eta < 1/r$}.
Our approach is to evaluate this integral directly.

We begin by expanding the contour to infinity, while avoiding the pole and cut
by wrapping around them in the usual way.
However, in order to deform the circular contour out to infinity, we must be able to continue
the integrand to the entire $z$-plane in a way that it remains analytic.
For any given
\mbox{$z\in\mathbb{C}$},
the solution of the corresponding quadratic equation has two distinct branches, but
it is {\it a priori} arbitrary how to assign the plus and minus signs, that multiply
the discriminant term, to the branches.
The most general form of $\zeta_\pm(z)$ is
\begin{equation}
\zeta_\pm(z) = \half\left[b(z) \pm \xi(z)\left(\Delta(z)\right)^{1/2}\right] \;,
\end{equation}
where $\xi(z)$ is an arbitrary function such that
\mbox{$\xi^2(z) \equiv 1$},
and the principal branch of the square-root function is adopted.
Analyticity of $\zeta_\pm(z)$ everywhere in the complex $z$-plane away from the
positive real axis is achieved by the choice
\begin{equation}
\xi(z) = \left\{
\begin{array}{ccc}
+1 & \quad\text{for}\quad & \Rez z \leq (x_+ + x_-)/2 \\
-1 & \quad\text{for}\quad & \Rez z > (x_+ + x_-)/2
\end{array}
\right. \;.
\label{xiz}
\end{equation}
With this definition, we observe that
\begin{equation}
\zeta_-(1/r) = \left\{
\begin{array}{ccc}
r               & \quad\text{for}\quad & r_{\text{hi}} > r^2 \\
r_{\text{hi}}/r & \quad\text{for}\quad & r_{\text{hi}} \leq r^2
\end{array}
\right. \;,
\end{equation}
which implies that there is no pole when
\mbox{$r_{\text{hi}} \geq r^2$}.
The behaviour of $\zeta_+(1/r)$ is the opposite and follows directly from (\ref{SumProd}).
With the choice in (\ref{xiz}), we also obtain the uniform asymptotic behaviour
\begin{equation}
\zeta_+(z) \asym{|z|\to\infty} -r_{\text{lo}}z \;, \quad
    \zeta_-(z) \asym{|z|\to\infty} -\frac{r_{\text{hi}}}{r_{\text{lo}}z} \;,
\end{equation}
which is sufficient to ensure the vanishing of the contour at infinity.

After deformation of the contour, $p_n$
may be expressed as the sum of a pole contribution and a cut integral
\mbox{$p_n = P_{\text{pol}}(n) + P_{\text{cut}}(n)$}, where
\begin{equation}
P_{\text{pol}}(n) = \Theta(r^2 - r_{\text{hi}})
     \left[1 - \frac{r(1-r)}{r_{\text{lo}}}\right]{\cdot}(1-r)r^{n-1} \;,
\label{MargPol}
\end{equation}
and
\begin{equation}
P_{\text{cut}}(n) = \frac{1-r}{2\pi r}\int_{x_-}^{x_+}\frac{dx}{x^{n+1}}\,
     \frac{\sqrt{(x_+-x)(x-x_-)}}{x - 1/r} \;.
\end{equation}
The change of integration variable $x\mapsto u$,
\mbox{$0 \leq u \leq 1$},
according to
\begin{equation}
x = x_- + (x_+ - x_-)u \;,
\end{equation}
yields
\begin{equation}
P_{\text{cut}}(n) = \frac{1-r}{2\pi r x_{\text{dif}}^n}\int_0^1 du\,
     \frac{\sqrt{u(1-u)}}{(u + a)^{n+1}(u + b)} \;,
\end{equation}
with
\begin{equation}
a \equiv x_-/x_{\text{dif}} \;, \quad b \equiv (x_- - 1/r)/x_{\text{dif}} \;,
     \quad x_{\text{dif}} \equiv x_+ - x_- \;,
\end{equation}
so that
\mbox{$1/x_{\text{dif}} = r(a-b)$}.
These parameters may also be conveniently expressed in terms of the model inputs as
\begin{equation}
a = \frac{1+r}{4\sqrt{r_{\text{hi}}}} - \frac{1}{2}\;, \quad
     b = \frac{1}{4}\left(\frac{r}{\sqrt{r_{\text{hi}}}} +
     \frac{\sqrt{r_{\text{hi}}}}{r}\right) - \frac{1}{2}\;.
\end{equation}

In terms of Appell's $F_1$ bivariate hypergeometric function, we may write
\begin{align}
\begin{aligned}
P_{\text{cut}}(n) &= \frac{1-r}{16rx_{\text{dif}}^n}{\cdot}\frac{1}{a^{n+1}b}{\cdot}
     \frac{1}{B(3/2,3/2)}\int_0^1 du\, \frac{\sqrt{u(1-u)}}{\left(1 + u/a\right)^{n+1}
     \left(1 + u/b\right)} \\
&= \frac{(1-r)r^{n-1}}{16ab}{\cdot}\left(1 - b/a\right)^n F_1(3/2,n+1,1,3;-1/a,-1/b) \;,
\end{aligned}
\end{align}
where it is useful to note that
\mbox{$1 - b/a = 1/(rx_-)$},
and with
\mbox{$B(\mu,\nu)$}
denoting the beta function.
This follows from Picard's integral representation \citep{NP:Picard81}
\begin{equation}
F_1(\lambda,\rho,\sigma,\lambda+\mu;u,v) = \frac{1}{B(\mu,\lambda)}
     \int_0^1dx\, \frac{x^{\lambda-1}(1-x)^{\mu-1}}
     {(1-ux)^\rho(1-vx)^\sigma} \;,
\end{equation}
as reproduced in \citep[sections 3.211, 9.18]{NP:Gradshteyn07}.
When, as in the present application, the arguments $u,v$ are real and negative,
the Appell $F_1$-hypergeometric function can be efficiently evaluated by means of a
Gauss-Jacobi quadrature rule \citep{NP:Gautschi04,NP:Gautschi05}.
Furthermore, in the present application, where
\mbox{$(\mu,\lambda) = (3/2,3/2)$},
the Gauss-Jacobi quadrature specializes to a Gauss-Chebyshev quadrature of the second kind.
To see this, we make the further change of integration variable
\mbox{$v = 2u - 1$}, which casts the integral into the form
\begin{equation}
P_{\text{cut}}(n) = C_n\int_{-1}^{+1} \frac{dv}{\pi/2}\, \frac{\sqrt{1 - v^2}}{(v+\alpha)^{n+1}(v+\beta)} \;,
     \quad C_n \equiv \frac{1-r}{4r}{\cdot}\left(\frac{2}{x_{\text{dif}}}\right)^n \;,
\label{pcutv}
\end{equation}
where we have introduced the parameters
\begin{equation}
\alpha \equiv 2a+1 \;, \quad \beta \equiv 2b + 1 \;.
\label{alphabeta}
\end{equation}
Thus,
\begin{equation}
\alpha = \frac{1+r}{2\sqrt{r_{\text{hi}}}} \;, \quad
\beta  = \frac{1}{2}\left(\frac{r}{\sqrt{r_{\text{hi}}}} +
     \frac{\sqrt{r_{\text{hi}}}}{r}\right) \;.
\end{equation}
It is known that Gauss-Chebyshev quadrature of the second kind
is equivalent to the trapezoidal rule on the unit circle \citep{NP:Chawla70}.
Accordingly, if we set
\mbox{$v = \cos(\pi\tau)$},
so that
\mbox{$0 \leq \tau \leq 1$},
and subsequently uniformly discretize the unit interval into $L$ sub-intervals
each terminated at the point
\mbox{$\tau_k = k/L$},
\mbox{$k = 1,2,\ldots,L$},
then we obtain the $L$-point quadrature rule for the cut integral given by
\begin{equation}
P_{\text{cut}}(n) \simeq C_n\sum_{k=1}^L \frac{w_k}{[\cos(\pi\tau_k)+\alpha]^{n+1}
     {\cdot}[\cos(\pi\tau_k)+\beta]} \;, \quad w_k \equiv (2/L){\cdot}\sin^2(\pi\tau_k) \;.
\end{equation}
This structure mirrors the discussion of the distribution of the partial busy period
for a multi-server Markovian queue found in \citep{NP:Zuk23A}.
An alternative way of expressing this, along the lines of \citep{NP:Zuk23A}, is
\begin{equation}
P_{\text{cut}}(n) \simeq \frac{x_{\text{dif}}}{4}{\cdot}\sum_{k=1}^L \frac{W_k}
     {\left[x_{\text{dif}}(U_k + a)\right]^{n+1}} \;,
\label{PGQ}
\end{equation}
with
\begin{equation}
U_k \equiv \cos^2(\pi\tau_k/2)\;, \quad
W_k \equiv \frac{2(1-r)}{rL}{\cdot}\frac{1 - U_k}{1 + b/U_k} \;.
\label{UWGQ}
\end{equation}
The case
\mbox{$\beta = 1$},
equivalently
\mbox{$b = 0$},
is problematic for the computation of the integral in (\ref{pcutv}) via the Gauss-Chebyshev
quadrature of the second kind as the denominator vanishes at an integration end-point,
where the weight function also vanishes.
Thus, a removable singularity appears.
At this critical point, which corresponds to parameter choice
\mbox{$r_{\text{hi}} = r^2$}
--- the same point at which the pole contribution disappears ---
(\ref{pcutv}) may be cast into a form suitable for Gauss-Chebyshev quadrature of the first kind:
\begin{equation}
P_{\text{cut}}(n) = C_n\int_{-1}^{+1} \frac{dv}{\pi/2}\, \frac{1}{\sqrt{1-v^2}}{\cdot}
     \frac{1 - v}{(v + \alpha)^{n+1}} \;.
\end{equation}
The loss of accuracy away from the critical point but close to it may be addressed by
noting that there is no impediment to evaluating (\ref{pcutv}) by Gauss quadrature of the
first kind in the general case, noting, in particular, that it does not place a node at
\mbox{$v = -1$}.
Furthermore, it is known that  Gauss-Chebyshev quadrature of the first kind is
equivalent to the mid-point rule on the unit circle.
In fact, the quadrature rule expressed by (\ref{PGQ}) and (\ref{UWGQ})
holds equally for the midpoint rule,
provided that one now chooses
\mbox{$\tau_k =  (k - 1/2)/L$}
for
\mbox{$k = 1,2,\ldots,L$}.
The advantage of the associations with the trapezoidal and mid-point rules is that
they can both be applied iteratively, in order to attain a pre-assigned level of
accuracy \citep{NP:Rice80}.
Generally, in a Gaussian quadrature scheme, the nodes and weights are generated for a
pre-selected quadrature order $L$, and it is difficult to estimate {\it a priori}
the accuracy that will result for any particular problem.
While one can successively increase the order to refine the integration grid,
thus achieving greater accuracy,
the nodes from the previous iteration cannot be re-used.
The mapping of a Gaussian quadrature to the trapezoidal rule allows node re-use if the
grid is doubled on each refinement and, for the mid-point rule, node re-use results
from tripling the grid on each refinement \citep{NP:Press07}.
While this makes the trapezoidal rule more efficient than the mid-point rule
in reaching a desired level of convergence, in the present problem,
a price is paid for the singular behaviour around the critical point.
It should be appreciated that the trapezoidal and mid-point rules do not generally
implement Gaussian quadratures.
The fact that it is so in the present case leads to an exponential rate of
convergence \citep{NP:Trefethen14}.

Finally, we may note that,
in the special case
\mbox{$r_{\text{hi}} = r^2$},
we have
\begin{equation}
p_n = \tensor[_2]{F}{_1}\left(n+1, 1/2; 2; -4r/(1-r)\right) {\cdot} r^n \;,
\end{equation}
recalling that the pole contribution vanishes here,
and using the result that the Gauss hypergeometric function has integral representation
\begin{equation}
\tensor[_2]{F}{_1}(a,b;c;z) = \frac{1}{B(b,c-b)}\int_0^1 dt\, \frac{t^{b-1}(1 - t)^{c-b-1}}
     {(1 - zt)^a} \;.
\end{equation}

\subsection{Asymptotic Limits}
Let us consider the integral
\begin{equation}
\mathcal{I}_n(a,b) \equiv \int_0^1 du\, \frac{\sqrt{u(1-u)}}{(u+a)^{n+1}(u+b)}
     \asym{n\to\infty} \frac{1}{b}\int_0^1 du\, \frac{\sqrt{u}}{(u+a)^{n+1}} \;,
\end{equation}
and continue the limiting procedure with
%
\begin{align}
\begin{aligned}
\int_0^1 du\, \frac{\sqrt{u}}{(u+a)^{n+1}} &\;\;=\;\; \frac{2}{a^{n-1/2}}\int_0^{a^2}
     dv\, \frac{v^2}{(1 + v^2)^{n+1}} \\
&\asym{n\to\infty} \frac{2}{a^{n-1/2}}\int_0^\infty
     dv\, \frac{v^2}{(1 + v^2)^{n+1}} \\
&\;\;=\;\; \frac{\sqrt{\pi}}{2a^{n-1/2}}{\cdot}\frac{\Gamma(n-1/2)}{\Gamma(n+1)} \;.
\end{aligned}
\end{align}
Hence, using the result that
\mbox{$\Gamma(n+c)\asym{n\to\infty} \Gamma(n)n^c$},
we obtain
\begin{equation}
\mathcal{I}_n(a,b) \asym{n\to\infty} \frac{\sqrt{\pi}}{2ba^{n-1/2}}{\cdot}\frac{1}{n^{3/2}} \;,
\end{equation}
and so we conclude that
\begin{equation}
P_{\text{cut}}(n) \asym{n\to\infty} \frac{1-r}{4rb}\sqrt{\frac{a}{\pi n^3}}{\cdot}
     \left[(1 - b/a)r\right]^n \;,
\end{equation}
provided
\mbox{$b\neq 0$}.
The situation where
\mbox{$b = 0$}
occurs when
\mbox{$r_{\text{hi}} = r^2$}.
In this case, we must consider
\begin{equation}
\mathcal{I}_n(a) \equiv \int_0^1 du\, \frac{1}{(u+a)^{n+1}}\sqrt{\frac{1-u}{u}}
     \asym{n\to\infty} \int_0^1 \frac{du}{\sqrt{u}}\, \frac{1}{(u+a)^{n+1}} \;.
\end{equation}
Now,
\begin{align}
\begin{aligned}
\int_0^1 \frac{du}{\sqrt{u}}\, \frac{1}{(u+a)^{n+1}} &\;\;=\;\; \frac{2}{a^{n+1/2}}\int_0^{a^2}
     dv\, \frac{1}{(1 + v^2)^{n+1}} \\
&\asym{n\to\infty} \frac{2}{a^{n+1/2}}\int_0^\infty
     dv\, \frac{1}{(1 + v^2)^{n+1}} \\
&\;\;=\;\; \frac{\sqrt{\pi}}{a^{n+1/2}}{\cdot}\frac{\Gamma(n+1/2)}{\Gamma(n+1)} \;,
\end{aligned}
\end{align}
and so it follows that
\begin{equation}
\mathcal{I}_n(a) \asym{n\to\infty} \frac{1}{a^{n+1/2}}\sqrt{\frac{\pi}{n}} \;.
\end{equation}
Therefore, when
\mbox{$r_{\text{hi}} = r^2$},
in which case
\mbox{$a = (1-r)/(4r)$},
we obtain
\begin{equation}
P_{\text{cut}}(n) \asym{n\to\infty} \sqrt{\frac{1-r}{\pi r}}{\cdot}\frac{r^n}{\sqrt{n}} \;.
\end{equation}

\section{Waiting-Time Distribution}
\label{WaitTime}
It is known that
the queue-length marginals can be recovered from the waiting-time marginals
by appealing to the distributional form of Little's law \citep{NP:Bertsimas95,NP:Keilson88},
bearing in mind the {\it caveat} that it is applicable to priority queues only within a given priority level.
Let us denote the Poisson PMF for arrival rate $\lambda$ and time $t$ by
\begin{equation}
\Lambda_n(\lambda t) \equiv \frac{(\lambda t)^n}{n!}e^{-\lambda t} \;,
\end{equation}
for arrival numbers
\mbox{$n = 0,1,2,\ldots$}.
Then, one has
\begin{equation}
P_\mathscr{L}(n) = \int_0^\infty dt\, P_\mathscr{W}(t)\Lambda_n(r_{\text{lo}}t) \;,
\label{DistLittle}
\end{equation}
where $\mathscr{L}, \mathscr{W}$ denote the queue-length and waiting-time
random variables (RVs), respectively.
This can be inverted, hence enabling us to derive the waiting-time probability density
function (PDF) $P_\mathscr{W}(t)$
as a by-product of our derivation of the queue-length marginal PMF $P_\mathscr{L}(n)$.
It is implicit in (\ref{DistLittle}) that we are dealing with the low-priority level
and that we have, without loss of generality, set
\mbox{$c\mu = 1$},
which means that
\mbox{$\lambda_{\text{lo}} = r_{\text{lo}}$}.
This simply implies that time is being measured in units of $c\mu$.

Let us consider
\begin{equation}
g_{\text{lo}}(z) = \sum_{n=0}^\infty z^n P_\mathscr{L}(n) = \int_0^\infty dt\, P_\mathscr{W}(t)
     K(z, r_{\text{lo}}t) \;,
\end{equation}
where $K(z,\lambda t)$
is the PGF of the Poisson distribution
\begin{equation}
K(z,\lambda t) \equiv \sum_{n=0}^\infty z^n\Lambda_n(\lambda t) = e^{\lambda t(z-1)}\;.
\end{equation}
It follows immediately that the moment generating function (MGF) of the waiting-time RV,
\begin{equation}
M_\mathscr{W}(z) \equiv \langle e^{z\mathscr{W}}\rangle_\mathscr{W} =
     \int_0^\infty dt\, P_\mathscr{W}(t) e^{zt} \;,
\end{equation}
is given by
\mbox{$M_\mathscr{W}(z) = g_{\text{lo}}(1 + z/r_{\text{lo}})$}.
The waiting-time PDF can now be recovered via the inverse Laplace transform as
\begin{equation}
P_\mathscr{W}(t) = \int_{\Gamma}\frac{dz}{2\pi i}\,
      e^{-tz} g_{\text{lo}}(1+z/r_{\text{lo}}) \;,
\end{equation}
with Bromwich contour $\Gamma$ such that
\mbox{$z = \gamma + iy$},
for
\mbox{$-\infty < y < +\infty$}
with any given real constant
$\gamma$ such that
\mbox{$0 < \gamma < (1/r-1)r_{\text{lo}}$}.
After a scaling and shift of the integration variable $z$, and a subsequent deformation
of the resulting Bromwich contour so that it wraps around the pole and cut of the function
\mbox{$g_{\text{lo}}(z)$}
on the real axis, we arrive at the representation
\begin{equation}
P_\mathscr{W}(t)dt = r_{\text{lo}}dt{\cdot}\int_{\mathcal{C}_{\text{pol}} \cup \mathcal{C}_{\text{cut}}}
     \frac{dz}{2\pi i}\, e^{-r_{\text{lo}}t(z-1)} g_{\text{lo}}(z) \;,
\label{WaitContour}
\end{equation}
noting that
\mbox{$\Rez z > 1$}
everywhere on
\mbox{$\mathcal{C} = \mathcal{C}_{\text{pol}} \cup \mathcal{C}_{\text{cut}}$}.
An alternative derivation is presented in Appendix~\ref{Laguerre}. While more complex,
it has the merit of establishing a novel representation for the waiting-time
PDF in terms of the queue-length PMF as a convergent series of Laguerre polynomials
whose coefficients are finite linear combinations of queue-length probabilities,
corresponding to a binomial transform.

Let us consider the closed contour integral, for an arbitrary analytic function
\mbox{$f(z)$},
given by
\begin{equation}
\mathcal{I}[f] \equiv \oint_\mathcal{C} \frac{dz}{2\pi i}\, f(z)g_{\text{lo}}(z) \;,
\end{equation}
with
\begin{equation}
g_{\text{lo}}(z) = \frac{1-r}{\zeta_+(z) - r} = \frac{1}{r_{\text{lo}}r}{\cdot}\frac{1-r}{z - 1/r}
     {\cdot}\left[\zeta_-(z) - r\right] \;,
\end{equation}
where we recall that
\begin{equation}
\lim_{z\to 1/r} \left(z - 1/r\right)g_{\text{lo}}(z) = \frac{1-r}{r^2}\left[\frac{r(1-r)}{r_{\text{lo}}} - 1\right]
     \Theta(r^2 - r_{\text{hi}}) \;.
\end{equation}
The chosen contour deformation leads to the pole/cut decomposition
\begin{equation}
\mathcal{I}[f] = \int_{\mathcal{C}_{\text{pol}} \cup \mathcal{C}_{\text{cut}}}
     \frac{dz}{2\pi i}\, f(z)g_{\text{lo}}(z)
     = \mathcal{I}_{\text{pol}}[f] + \mathcal{I}_{\text{cut}}[f] \;,
\end{equation}
with
\begin{equation}
\mathcal{I}_{\text{pol}}[f] = \frac{1-r}{r^2}\left[1 - \frac{r(1-r)}{r_{\text{lo}}}\right]
     \Theta(r^2 - r_{\text{hi}}){\cdot}f(1/r) \;,
\end{equation}
and
\begin{align}
\begin{aligned}
\mathcal{I}_{\text{cut}}[f] &= \frac{1-r}{2\pi r}\int_{x_-}^{x_+} dx\,
     \frac{\sqrt{(x_+ - x)(x - x_-)}}{x - 1/r} f(x) \\
&= \frac{(1-r)x_{\text{dif}}}{2\pi r}\int_0^1 du\,\frac{\sqrt{u(1-u)}}{u+b}
     f(x_{\text{dif}}(u+a)) \;.
\end{aligned}
\end{align}
It should be noted that $\mathcal{C}_{\text{pol}}$ and $\mathcal{C}_{\text{cut}}$
are both clockwise contours about the pole at
\mbox{$z = 1/r$}
and cut along
\mbox{$z \in [x_-,x_+]$},
respectively.
Application to the waiting-time PDF entails considering
\begin{equation}
f(z;t) = r_{\text{lo}} e^{-r_{\text{lo}}t(z-1)} \;,
\end{equation}
in which case
the wait-conditional waiting-time PDF, given by
\mbox{$P_W(t) = \mathcal{I}[f](t)$},
becomes
\begin{align}
\begin{aligned}
P_\mathscr{W}(t) &= -(1-r)\left(\frac{r_{\text{hi}}}{r^2}-1\right)\Theta(r^2-r_{\text{hi}})
     e^{-tr_{\text{lo}}(1/r-1)} \\
     &\quad {}+ \frac{2(1-r)\sqrt{r_{\text{hi}}}}{\pi r}e^{-t(4\sqrt{r_{\text{hi}}}a - r_{\text{lo}})}
     \int_0^1 du\, e^{-4t\sqrt{r_{\text{hi}}}{\cdot}u}{\cdot}\frac{\sqrt{u(1-u)}}{u+b} \;,
\end{aligned}
\label{WaitPDF}
\end{align}
where we have observed that
\mbox{$r_{\text{lo}}x_{\text{dif}} = 4\sqrt{r_{\text{hi}}}$}.
The integral appearing here is an Euler-type integral representation for a special case
of the Humbert series $\Phi_1$, given by \cite{NP:Erdelyi40} as
\begin{equation}
\Phi_1(a,b,c;x,y) = \frac{1}{B(a,c-a)}\int_0^1 dt\,e^{yt} {\cdot}\frac{t^{a-1}(1-t)^{c-a-1}}
     {(1 - xt)^b} \;.
\end{equation}
The Humbert series, also known as confluent Appell functions,
are bivariate generalizations of Kummer's confluent hypergeometric function
$\tensor[_1]{F}{_1}(a;c;y)$.
In the present case, we have
\mbox{$\tensor[_1]{F}{_1}(a;c;y) = \Phi_1(a,b,c;0,y)$}
for any $b$.
The function $\Phi_1$ in the general case is efficiently computed via Gauss-Jacobi
quadrature \citep{NP:Gautschi02}.
In the present specific application, the iterative Gauss-Chebyshev scheme described
earlier can be utilized effectively.

It is straightforward to show explicitly that
\mbox{$P_\mathscr{W}(0) = 1-r$}.
A useful intermediate result is
\begin{equation}
\left(\sqrt{b+1} - \sqrt{b}\right)^2 = \left\{
\begin{array}{ccc}
\sqrt{r_{\text{hi}}}/r & \quad\text{for}\quad  & r_{\text{hi}} \leq r^2 \;, \\
r/\sqrt{r_{\text{hi}}} & \quad\text{for}\quad  & r_{\text{hi}} >    r^2 \;.
\end{array}
\right.
\end{equation}
By contrast, it requires considerable algebra to confirm explicitly that
\mbox{$P_\mathscr{W}(t)$}
integrates to unity; whereas it is easy to see this on general grounds.
Taking the $t$-integral inside the $z$-integral yields
\begin{equation}
\int_0^\infty dt\, f(z;t) = \frac{1}{z-1} \;,
\end{equation}
provided
\mbox{$\Rez z > 1$}.
Then, on deforming the integration contour
\mbox{$\mathcal{C}\mapsto \mathcal{C}_\eta$},
with the contour $\mathcal{C}_\eta$ denoting the anti-clockwise circle
around the origin of radius
\mbox{$1 < \eta< 1/r$},
we obtain
\begin{equation}
\int_0^\infty dt\, P_\mathscr{W}(t) = \oint_{\mathcal{C}_\eta} \frac{dz}{2\pi i}\,
     \frac{g_{\text{lo}}(z)}{z-1} = g_{\text{lo}}(1) = 1 \;.
\end{equation}

The wait-conditional low-priority waiting-time marginal PDF given by (\ref{WaitPDF})
agrees with the result in \citep{NP:Davis66} derived from a different starting point,
based on a first-principles derivation of the MGF.
An alternative derivation of the same MGF has been presented in \citep{NP:Kella85}.
This result can be used to represent waiting-time marginals for all priority levels
in the general multi-level problem.
Let $r_\kappa$,
\mbox{$\kappa = 1,2,\ldots,K$}
denote the level traffic intensities for a $K$-level problem,
ordered in descending priority, so that
\mbox{$r_{\text{agg}} = r_1 + r_2 + \cdots + r_K < 1$}
is the total (aggregated) traffic intensity.
A common mean service rate $\mu$ is assumed.
We also introduce the cumulative quantities \citep{NP:Holley54}
\begin{equation}
\sigma_\kappa \equiv \sum_{j=1}^\kappa r_j \;,
\end{equation}
so that
\mbox{$\sigma_K = r_{\text{agg}}$},
and make the associations with an effective two-level problem
according to
\begin{equation}
r_{\text{lo}} = r_\kappa \;, \quad r_{\text{hi}} = \sigma_{\kappa-1} \;, \quad
     r = r_{\text{lo}} + r_{\text{hi}} \;,
\end{equation}
Then $P_{\mathscr{W}_\kappa}(t)$,
as given by (\ref{WaitPDF}),
yields the wait-conditional waiting-time PDF for the $\kappa$-th highest
priority level.
In particular, for the highest priority level
(\mbox{$\kappa = 1$}),
we would set
\mbox{$r_{\text{lo}} = r_1$},
\mbox{$r_{\text{hi}} = 0$},
so that
\mbox{$r = r_1$},
which leads to the result, arising solely from the pole contribution, that
\begin{equation}
P_{\mathscr{W}_1}(t) = (1 - r_1)e^{-(1-r_1)t} \;,
\end{equation}
being an exponential distribution, as expected.
Since the waiting-time marginal PDFs have already appeared in the literature
in a similar form \citep{NP:Davis66},
we shall not discuss them further in the present work.
As far as we are aware, the only explicit representation of the waiting-time
distributions that exist in the literature is the integral form
first derived in \citep{NP:Davis66}
and elaborated upon here. In Appendix~\ref{Laguerre}, we establish a different representation
as a Laguerre series.

\section{Joint Distribution}
\label{Joint}
For the wait-conditional joint queue-length distribution,
the pole/cut decomposition reads
\begin{equation}
P(\ell,m) = P_{\text{pol}}(\ell,m) + P_{\text{cut}}(\ell,m) \;,
\end{equation}
with
\begin{align}
\begin{aligned}
P_{\text{pol}}(\ell,m) &= -\lim_{z\to 1/r}\left(z - 1/r\right)g_\ell(z)/z^{m+1} \;, \\
P_{\text{cut}}(\ell,m) &= \frac{1}{2\pi i}\int_{x_-}^{x_+}\frac{dx}{x^{m+1}}\,
     \left[g_\ell(x+i\epsilon) - g_\ell(x-i\epsilon)\right] \;,
\end{aligned}
\end{align}
where $\epsilon$ is a positive infinitesimal.
On inserting the explicit form of the PGF, the cut contribution reads
\begin{equation}
P_{\text{cut}}(\ell,m) = \frac{1-r}{2\pi i}\int_{x_-}^{x_+} dx\,
     \frac{i\sqrt{|\Delta(x)|}}{x^{m+1}(1 - rx)} \chi_\ell(x)
     - \frac{1-r}{2\pi i}\int_{x_-}^{x_+} dx\,
     \frac{i\sqrt{|\Delta(x)|}}{x^{m}(1 - rx)} \chi_{\ell+1}(x) \;,
\end{equation}
where
\begin{equation}
\chi_\ell(x) \equiv -\frac{\zeta^{\prime\ell}_-(x) - \zeta^{\prime\ell}_+(x)}
     {i\sqrt{|\Delta(x)|}} \;, \quad
     \zeta^{\prime}_\pm(x) \equiv \lim_{\epsilon\to 0^+} \zeta_\pm(x-i\epsilon) \;,
\end{equation}
is an analytic function.
On the cut, we can write
\begin{align}
\begin{aligned}
\zeta_\pm(x) &= \frac{r_{\text{lo}}}{4}\left(\sqrt{x_+-x} \pm i\sqrt{x-x_-}\right)^2 \\
     &= \sqrt{r_{\text{hi}}}\left(\sqrt{1-u} \pm i\sqrt{u}\right)^2 \;.
\end{aligned}
\end{align}
We also have
\begin{equation}
|\Delta(x)| = r_{\text{lo}}^2 x_{\text{dif}}^2 u(1-u)
     = 16r_{\text{hi}}u(1-u) \;.
\end{equation}
This leads to the result
\begin{align}
\begin{aligned}
\chi_\ell(x(u)) &= (-1)^{\ell-1}r_{\text{hi}}^{(\ell-1)/2}{\cdot}
     \frac{\left(y + \sqrt{y^2-1}\right)^\ell - \left(y - \sqrt{y^2-1}\right)^\ell}
     {2\sqrt{y^2-1}} \\
     &= (-1)^{\ell-1}r_{\text{hi}}^{(\ell-1)/2}{\cdot}U_{\ell-1}(y) \\
     &= (-1)^{\ell-1}r_{\text{hi}}^{(\ell-1)/2}{\cdot}U_{\ell-1}(2u-1) \\
     &= r_{\text{hi}}^{(\ell-1)/2}{\cdot}\frac{\sin(2\ell\asin\sqrt{u})}
     {\sin(2\asin\sqrt{u})} \;,
\end{aligned}
\end{align}
where we had temporarily set
\mbox{$y = 2u-1$}
and used some standard representations for the Chebyshev polynomial $U_\ell(y)$ of the
second kind of order $\ell$.
One may observe that
\mbox{$\chi_0(x) = 0$},
\mbox{$\chi_1(x) = 1$}.
From the Chebyshev generating function
\begin{equation}
\sum_{\ell=0}^\infty t^\ell U_\ell(y) = \frac{1}{1 - 2ty + t^2} \;,
\end{equation}
one obtains
\begin{equation}
\sum_{\ell=0}^\infty \chi_\ell(x(u)) = \frac{1}{1 + 2\sqrt{r_{\text{hi}}}(2u-1) + r_\text{hi}} \;,
\end{equation}
which is useful for verifying that the low-priority marginal is recovered
when one sums over the high-priority indices $\ell$.

After some algebraic manipulation, we arrive at the expression
\begin{equation}
\begin{split}
P_{\text{cut}}(\ell,m) = &-C_1\int_0^1\frac{du}{\pi}{\cdot}\frac{\sqrt{u(1-u)}}{(u+a)^{m+1}(u+b)}{\cdot}
     \frac{\sin(2\ell\asin\sqrt{u})}{\sin(2\asin\sqrt{u})} \\
     &{}+ C_2\int_0^1\frac{du}{\pi}{\cdot}\frac{\sqrt{u(1-u)}}{(u+a)^{m}(u+b)}{\cdot}
     \frac{\sin(2(\ell+1)\asin\sqrt{u})}{\sin(2\asin\sqrt{u})} \;,
\end{split}
\end{equation}
with coefficients
\begin{equation}
C_1 \equiv \frac{2(1-r)}{r}\left(\frac{r_{\text{lo}}}{4}\right)^{m+1}
     {\cdot}r_{\text{hi}}^{(\ell-m-1)/2} \;, \quad
C_2 \equiv \frac{2(1-r)}{r}\left(\frac{r_{\text{lo}}}{4}\right)^{m}
     {\cdot}r_{\text{hi}}^{(\ell-m+1)/2} \;.
\end{equation}
As for the marginal PMF, we set
\mbox{$u = (1+v)/2$},
\mbox{$v  = \cos(\pi\tau)$},
\mbox{$0 \leq \tau \leq 1$},
so that
\begin{equation}
u(\tau) = \left[1 + \cos(\pi\tau)\right]/2 = \cos^2(\pi\tau/2) \;,
\end{equation}
and
\begin{equation}
\frac{du}{\pi}{\cdot}\frac{\sqrt{u(1-u)}}{u} = -d\tau{\cdot}[1 - u(\tau)] \;.
\end{equation}
We also have
\begin{equation}
C_\ell(u) \equiv \frac{\sin(2\ell\asin\sqrt{u})}{\sin(2\asin\sqrt{u})}
     = (-1)^{\ell-1}\frac{\sin(\ell\pi\tau)}{\sin(\pi\tau)}
     = (-1)^{\ell-1}\ell{\cdot}\frac{\sinc(\ell\pi\tau)}{\sinc(\pi\tau)} \;,
\end{equation}
and we can write
\mbox{$P_{\text{cut}}(\ell,m) = \Lambda(\ell+1,m) - \Lambda(\ell,m+1)$},
where
\begin{equation}
\Lambda(\ell,m) \equiv \frac{2(1-r)}{r}{\cdot}\frac{r_{\text{hi}}^{\ell/2}}{x_-^m}
     \int_0^1 d\tau\, \frac{[1 - u(\tau)]C_\ell(u(\tau))}{[1 + u(\tau)/a]^m[1 + b/u(\tau)]} \;.
\end{equation}
The $L$-point quadrature with nodes $\tau_k$,
\mbox{$k = 1,2,\ldots,L$},
corresponding to the trapezoidal rule is obtained by introducing
the integration grid
\mbox{$\tau_k = k/L$},
while the mid-point rule corresponds to
\mbox{$\tau_k = (k-1/2)/L$}.
In both cases, the resulting quadrature rule may be expressed as
\begin{equation}
\Lambda(\ell,m) \simeq r_{\text{hi}}^{\ell/2}\sum_{k=1}^L\frac{W_kC_\ell(U_k)}{[x_-(1 + U_k/a)]^m} \;,
\label{LamQuad}
\end{equation}
where the quadrature nodes and weights $U_k$, $W_k$, respectively, are given by (\ref{UWGQ}).
As a sanity check, we may note that the identity for the exclusively-low probabilities
\mbox{$P_{\text{cut}}(0,m+1) = r_{\text{lo}}P_{\text{cut}}(m)$}
implies that
\mbox{$\Lambda(1,m+1) = r_{\text{lo}}P_{\text{cut}}(m)$}.
Recalling that
\mbox{$x_- = ax_{\text{dif}}$},
it is straightforward to verify that the quadrature rules (\ref{LamQuad}) and (\ref{PGQ})
applied to the left-hand and right-hand sides, respectively,
reproduce this equality.
The explicit result for the pole contribution is easily obtained as
\begin{equation}
P_{\text{pol}}(\ell,m) = \Theta(r^2 - r_{\text{hi}}){\cdot}
     (1-r)r^m\left(1 - \frac{r_{\text{hi}}}{r^2}\right)
     \left(\frac{r_{\text{hi}}}{r}\right)^\ell \;.
\label{pjointpol}
\end{equation}
%

\section{Exact Expressions}
\label{Exact}
\subsection{Marginal Distribution}
In order to facilitate an exact representation of the low-priority marginal PMF,
we introduce a family of integrals to which we shall refer as the
marginal integrals. These are given by
\begin{equation}
M_n(\alpha,\beta) \equiv \int_0^\pi d\theta\, \frac{1 - \cos^2\theta}
     {(\alpha+\cos\theta)^{n+1}(\beta+\cos\theta)} \;,
\end{equation}
for
\mbox{$n = -1,0,\ldots$},
so that
\mbox{$M_{-1}(\alpha,\beta) = \pi(\beta - \sqrt{\beta^2 - 1})$}.
It is assumed that
\mbox{$|\alpha|, |\beta| > 1$}.
In terms of the marginal integral, we have
\begin{equation}
P_{\text{cut}}(n) = \frac{(1-r)r^{n-1}}{2\pi}{\cdot}
     (\alpha-\beta)^n M_n(\alpha,\beta) \;,
\end{equation}
with the parameters
\mbox{$\alpha,\beta$}
defined in (\ref{alphabeta}).
The marginal integral may be evaluated in terms of Legendre polynomials
\mbox{$P_n(\chi)$}
due to the identity \citep[section 3.661.4]{NP:Gradshteyn07}
\begin{equation}
\frac{1}{\pi}\int_0^\pi \frac{d\theta}{(\alpha + \cos\theta)^{n+1}} =
     \frac{1}{(\alpha^2 - 1)^{(n+1)/2}}P_n\left(\frac{\alpha}{\sqrt{\alpha^2 - 1}}\right) \;.
\end{equation}
for
\mbox{$\alpha > 1$}
and
\mbox{$n = -1,0,1,\ldots$},
where we formally set
\mbox{$P_{-1}(\chi) \equiv 1$}.
After a partial fraction decomposition of the integrand, one obtains
for
\mbox{$n = 0,1,\ldots$},
\begin{equation}
M_n(\alpha,\beta) = \frac{\pi}{(\alpha-\beta)^{n+1}}\sum_{k=-1}^n \psi_k(n)
     P_k\left(\frac{\alpha}{\sqrt{\alpha^2 - 1}}\right) \;,
\label{MP}
\end{equation}
with
\begin{equation}
\psi_k(n) = \frac{(\alpha - \beta)^k}{(\alpha^2-1)^{(k+1)/2}}
     \left[\beta^2-1 + (\alpha^2-\beta^2)\delta_{k,n} - (\alpha-\beta)^2\delta_{k,n-1}\right] \;,
\end{equation}
for
\mbox{$k = 1,2,\ldots,n$},
and
\begin{align}
\begin{aligned}
\psi_{-1}(n) &= (\beta-\alpha)\delta_{n,0} \;, \\
\psi_0(n)    &= (\beta^2-1)\left(\frac{1}{\sqrt{\alpha^2-1}} - \frac{1}{\sqrt{\beta^2-1}}\right)
     + \frac{\alpha^2-\beta^2}{\sqrt{\alpha^2-1}}\delta_{n,0}
     - \frac{(\alpha-\beta)^2}{\sqrt{\alpha^2-1}}\delta_{n,1}\;.
\end{aligned}
\end{align}
An immediate consequence is that
\begin{equation}
M_0(\alpha,\beta)/\pi = \frac{\psi_{-1}(0) + \psi_0(0)}{\alpha - \beta}
     = \frac{\alpha+\beta}{\sqrt{\alpha^2-1} +\sqrt{\beta^2-1}} - 1 \;,
\end{equation}
which, as expected, is symmetric in $\alpha$ and $\beta$.

The representation of the marginal PMF arising from (\ref{MP}) is clearly singular in the limit
\mbox{$\alpha - \beta \to 0^+$},
noting that, in the present queueing application, we always have
\mbox{$\alpha > \beta > 1$}.
This limit is encountered as
\mbox{$r_{\text{lo}}\to 0$}
or, equivalently,
\mbox{$\nu\to 1$}.
It is also singular in the limit
\mbox{$\alpha\to 1^+$},
which is encountered when
\mbox{$r_{\text{lo}}\to 0$},
\mbox{$r_{\text{hi}}\to 1$}
simultaneously.
The representation may be re-arranged and scaled to make it numerically robust
throughout the entire ergodic region.
The final algorithm is detailed in the implementation summary at the end of this section.

\subsection{Joint Distribution}
In order to facilitate an exact representation of the joint marginal PMF,
we introduce a family of integrals to which we shall refer as the
joint integrals. These are defined by
\begin{align}
\begin{aligned}
J_m^\ell(\alpha,\beta) &\equiv \int_0^\pi d\theta\, \frac{1 - \cos^2\theta}
     {(\alpha+\cos\theta)^{m+1}(\beta+\cos\theta)}U_\ell(\cos\theta) \\
&= \int_0^\pi d\theta\, \frac{\sin\theta{\cdot}\sin((\ell+1)\theta)}
     {(\alpha+\cos\theta)^{m+1}(\beta+\cos\theta)} \;.
\end{aligned}
\label{JIntg}
\end{align}
A relationship with the marginal integrals is given by
\begin{align}
\begin{aligned}
J_m^0(\alpha,\beta) &= M_m(\alpha,\beta) \;, \\
J_m^1(\alpha,\beta) &= 2M_{m-1}(\alpha,\beta) - 2\alpha M_m(\alpha,\beta) \;.
\end{aligned}
\end{align}
In the general case, the recurrence relation for the Chebyshev polynomial
\mbox{$U_\ell(\cos\theta)$}
implies the recurrence relation for the joint integral
\begin{equation}
J_m^{\ell+1}(\alpha,\beta) = 2J_{m-1}^\ell(\alpha,\beta)
     - 2\alpha J_m^\ell(\alpha,\beta) - J_m^{\ell-1}(\alpha,\beta) \;,
\end{equation}
for
\mbox{$\ell = 0,1,\ldots$},
where we take the seed values
\begin{equation}
J_m^{-1}(\alpha,\beta) \equiv 0 \;, \quad
     J_m^0(\alpha,\beta) = M_m(\alpha,\beta) \;,
\end{equation}
for
\mbox{$m = -1,0,1,\ldots$},
noting that
\mbox{$M_{-1}(\alpha,\beta) = \pi(\beta - \sqrt{\beta^2-1})$}.

Next, we introduce the intermediate integrals
\begin{align}
\begin{aligned}
I_m^\ell(\alpha,\beta) &\equiv \int_0^\pi d\theta\, \frac{\cos(\ell\theta)}
     {(\alpha + \cos\theta)^{m+1}(\beta + \cos\theta)} \;, \\
I_m^\ell(\alpha) &\equiv \int_0^\pi d\theta\, \frac{\cos(\ell\theta)}
     {(\alpha + \cos\theta)^{m+1}} \;,
\end{aligned}
\end{align}
for
\mbox{$m = -1,0,\ldots$},
\mbox{$\ell\in\mathbb{Z}$}.
One may observe that
\mbox{$M_n(\alpha,\beta) = \left[I_n^0(\alpha,\beta) - I_n^2(\alpha,\beta)\right]/2$}
and that
\mbox{$I_{-1}^\ell(\alpha,\beta) = I_0^\ell(\beta)$},
\mbox{$I_{-1}^\ell(\alpha) = \pi\delta_{\ell 0}$}.
A partial fraction decomposition yields
\begin{equation}
I_m^\ell(\alpha,\beta) = c_m I_0^\ell(\beta) - \sum_{k=0}^m c_k I_{m-k} ^\ell(\alpha) \;,
\end{equation}
where we have set
\mbox{$c_m \equiv 1/(\alpha - \beta)^{m+1}$}.
The joint integral may be expressed in terms of the intermediate integrals as
\begin{equation}
J_m^\ell(\alpha,\beta) =  I_m^\ell(\alpha,\beta) + \beta I_m^{\ell+1}(\alpha,\beta)
     - I_m^{\ell+1}(\alpha) \;.
\end{equation}
In terms of the joint integrals, we can write
\begin{equation}
P_{\text{cut}}(\ell, m) = \frac{(1-r)r_{\text{lo}}}{2\pi r}\left[
     \frac{(-\sqrt{r_{\text{hi}}})^{\ell}}{(x_{\text{dif}}/2)^{m-1}}
     J_{m-1}^{\ell}(\alpha,\beta) -
     \frac{(-\sqrt{r_{\text{hi}}})^{\ell-1}}{(x_{\text{dif}}/2)^m}
     J_{m}^{\ell-1}(\alpha,\beta)
     \right] \;.
\end{equation}

\subsection{$Q$-Functions}
To proceed further, we define a family of functions
\mbox{$Q_m^\ell(\chi)$},
that generalize the Legendre polynomials
and are related closely to the associated Legendre functions,
by means of the integral representation
\begin{equation}
Q_m^\ell(\chi) \equiv \frac{(-1)^\ell}{\pi}\left(\frac{\chi+1}{\chi-1}\right)^{\ell/2}
     \int_0^\pi d\theta\, \frac{\cos(\ell\theta)}{(\chi + \sqrt{\chi^2-1}\cos\theta)^{m+1}} \;.
\end{equation}
This may be compared with the integral representation of the associated Legendre function
given by \citep[section 8.711.2]{NP:Gradshteyn07}
\begin{equation}
P_m^\ell(\chi) = (-1)^\ell\frac{\Gamma(m+1)}{\Gamma(m-\ell + 1)}
     \int_0^\pi\frac{d\theta}{\pi}\,
     \frac{\cos(\ell\theta)}{(\chi + \sqrt{\chi^2 - 1}\cos\theta)^{m+1}} \;,
\end{equation}
noting also that the Legendre polynomials are recovered according to
\mbox{$P_m(\chi) = P_m^0(\chi)$}.
Clearly
\mbox{$Q_{-1}^{\ell}(\chi) = \delta_{\ell 0}$},
\mbox{$Q_{0}^{\ell}(\chi) = 1$}
and
\mbox{$Q_m^0(\chi) = P_m(\chi)$}.
It also follows that
\begin{equation}
Q_m^1(\chi) = (\alpha + \sqrt{\alpha^2-1})\left[\alpha P_m(\chi) -
     \sqrt{\alpha^2-1} P_{m-1}(\chi)\right] \;,
\end{equation}
for
\mbox{$m = 0,1,\ldots$}.
From their definition and the trigonometric identity
\begin{equation}
\cos[(\ell+1)\theta] + \cos[(\ell-1)\theta] = 2\cos(\ell\theta)\cos\theta \;,
\end{equation}
it is easily shown
that the $Q$-functions satisfy the recurrence relation
\begin{equation}
Q_m^{\ell+1}(\chi) = \frac{2\chi}{\chi-1}Q_m^{\ell}(\chi)
     - \frac{2}{\chi-1}Q_{m-1}^{\ell}(\chi)
     - \frac{\chi+1}{\chi-1}Q_m^{\ell-1}(\chi) \;,
\end{equation}
for
\mbox{$m = 0,1,\ldots$},
\mbox{$\ell\in\mathbb{Z}$}.
An explicit expansion of the $Q$-function is given by
\begin{equation}
Q_m^\ell(\chi) = \sum_{k=0}^m\binom{m+\ell}{k+\ell}\binom{k+m}{m}
     \left(\frac{\chi-1}{2}\right)^k \;.
\end{equation}
The relationship with associated Legendre functions is established by noting that,
for
\mbox{$\ell \leq m$},
\begin{equation}
Q_m^\ell(\chi) = \frac{(m-\ell)!}{m!}\left(\frac{\chi+1}{\chi-1}\right)^{\ell/2}
     P_m^\ell(\chi) \;.
\end{equation}
Another recurrence relation,
which follows directly from the recurrence relations for the associated Legendre functions, is
given by
\begin{equation}
Q_{m+1}^\ell(\chi) = \frac{2m+1}{m+1}\chi Q_m^\ell(\chi)
     - \frac{m^2-\ell^2}{m(m+1)}Q_{m-1}^\ell(\chi) \;,
\end{equation}
for
\mbox{$m = 1,2,\ldots$},
with the first few functions being
\begin{equation}
Q_0^\ell(\chi) = 1 \;, \quad
Q_1^\ell(\chi) = \chi + \ell \;, \quad
Q_2^\ell(\chi) = \tfrac{3}{2}\chi^2 + \tfrac{3}{2}\ell\chi + \tfrac{1}{2}(\ell^2-1) \;.
\end{equation}
This recurrence is particularly useful for practical purposes since, while superficially appearing
to be a two-dimensional recurrence, it is actually a family of independent one-dimensional recurrences
in $m$, each indexed by the value of $\ell$.

The generating function for the $Q$-functions is found to be
\begin{equation}
\sum_{m=0}^\infty z^m Q_m^\ell(\chi) = \frac{1}{\sqrt{z^2 - 2\chi z + 1}}
     \left(\frac{\chi - z - \sqrt{z^2 - 2\chi z + 1}}{\chi - 1}\right)^\ell \;.
\end{equation}
Using
\begin{equation}
\chi = \frac{\alpha}{\sqrt{\alpha^2-1}} \;, \quad
     z = \frac{\alpha - \beta}{\sqrt{\alpha^2-1}} \;,
\end{equation}
and noting that
\mbox{$z^2 - 2\chi z + 1 = (\beta^2-1)/(\alpha^2-1)$},
this may be cast in the alternative form
\begin{equation}
\sum_{m=0}^\infty z^m Q_m^\ell(\chi) = \sqrt{\frac{\alpha^2-1}{\beta^2-1}}{\cdot}
     \frac{d_\ell(\beta)}{d_\ell(\alpha)} \;,
\end{equation}
where we have set
\begin{equation}
d_\ell(\alpha) \equiv (\alpha - \sqrt{\alpha^2-1})^\ell = \frac{1}{(\alpha + \sqrt{\alpha^2-1})^\ell} \;,
\end{equation}
and analogously for
\mbox{$d_\ell(\beta)$}.

The single-argument intermediate integral can be expressed in terms of the $Q$-function as
\begin{align}
\begin{aligned}
I_m^\ell(\alpha) &= \frac{\pi(-1)^\ell}{(\alpha + \sqrt{\alpha^2-1})^\ell(\alpha^2-1)^{(m+1)/2}}
     Q_m^\ell\left(\frac{\alpha}{\sqrt{\alpha^2-1}}\right) \\
&= \pi(-1)^\ell\left(\frac{\chi-1}{\chi+1}\right)^{\ell/2}{\cdot}(\chi^2-1)^{(m+1)/2}
     Q_m^\ell(\chi) \;.
\end{aligned}
\end{align}
Consequently, we may express the joint integral in terms of the $Q$-functions,
according to
\begin{equation}
\begin{split}
\frac{(-1)^\ell J^\ell_m(\alpha,\beta)}{\pi c_m} &= d_{\ell+1}(\beta) + d_{\ell+1}(\alpha)
     \left(\frac{\alpha-\beta}{\sqrt{\alpha^2-1}}\right)^{m+1}Q_m^{\ell+1}(\chi) \\
     &\quad {}- \frac{1}{\sqrt{\alpha^2-1}}\sum_{k=0}^m
     \left(\frac{\alpha-\beta}{\sqrt{\alpha^2-1}}\right)^k\left[
     d_\ell(\alpha)Q_k^\ell(\chi) - \beta d_{\ell+1}(\alpha)Q_k^{\ell+1}(\chi)\right] \;.
\end{split}
\end{equation}

\subsection{Implementation}
We conclude this section by giving complete, self-contained
specifications of the algorithms
for obtaining exact closed-form algebraic expressions for the marginal and joint
queue-length PMFs that are also numerically efficient and robust.
Let us first recall that
\begin{equation}
\alpha = \frac{r+1}{2\sqrt{r_{\text{hi}}}} \;, \quad
     \beta = \frac{1}{2}\left(\frac{r}{\sqrt{r_{\text{hi}}}} + \frac{\sqrt{r_{\text{hi}}}}{r}\right) \;,
\end{equation}
and that
\begin{equation}
\chi = \frac{\alpha}{\sqrt{\alpha^2 - 1}} \;, \quad
     z = \frac{\alpha - \beta}{\sqrt{\alpha^2 - 1}} \;,
\end{equation}
so that
\begin{equation}
\chi = \frac{1+r}{\sqrt{(1+r)^2 - 4r_{\text{hi}}}} \;, \quad
     z = \frac{1 - r_{\text{hi}}/r}{\sqrt{(1+r)^2 - 4r_{\text{hi}}}} \;, \quad
     t \equiv \chi z \;.
\end{equation}
It is useful to note that
\begin{equation}
d_1(\alpha) = \frac{1}{2\sqrt{r_{\text{hi}}}}\left[1 + r - \sqrt{(1 + r)^2 - 4r_{\text{hi}}}\right]
     = \zeta_-(0)/\sqrt{r_{\text{hi}}} \;,
\label{d1alpha}
\end{equation}
and that
\begin{equation}
d_1(\beta) = \min\left(r/\sqrt{r_{\text{hi}}}, \sqrt{r_{\text{hi}}}/r\right) \;.
\label{d1beta}
\end{equation}
It is also useful to note the bounds on $z$ given by
\begin{equation}
0 \leq \frac{\alpha-\beta}{\sqrt{\alpha^2-1}} \leq \frac{1}{1+r} < 1 \;,
\end{equation}
which shows that, while $z$ never becomes large, it can be very small.
The lower and upper endpoints for $z$ correspond to
\mbox{$r_{\text{lo}} = 0$}
and
\mbox{$r_{\text{hi}} = 0$},
respectively.
We also have
\begin{align}
\begin{aligned}
1 &\leq \chi \leq (1+r)/(1-r) \;, \\
0 &\leq \chi z \leq 1/(1+r) \;.
\end{aligned}
\end{align}
Table~\ref{tab:sing} enumerates the singular limits of the marginal and
joint integrals.
The algorithms described in the following two subsections implement scalings and
algebraic manipulations to mitigate their effects, and render the final expression for the
marginal and joint PMF numerically robust across the entire ergodic region.

\begin{table}
\TABLE
{Singular Regions\label{tab:sing}}
{\begin{tabular}{|c|c|c|}
\hline
$\alpha,\beta$ & $r_{\text{lo}},r_{\text{hi}}$ & $z$ \\
\hline\hline
$\beta\to 1^+$             & $r_{\text{hi}}\to r^2$                   &                \\
$\alpha - \beta \to 0^+$   & $r_{\text{lo}}\to 0$                     & $z\to 0$       \\
$\alpha - \beta \to\infty$ & $r_{\text{hi}}\to 0$                     & $z\to 1/(1+r)$ \\
$\alpha\to 1^+$            & $r_{\text{lo}}\to 0, r_{\text{hi}}\to 1$ & $z\to 0$       \\
\hline
\end{tabular}}
{}
\end{table}

\subsubsection{Marginal Distribution}
\label{MargDistLeg}
We introduce scaled Legendre polynomials
\mbox{$\hat{P}_n(t;z) \equiv z^n P_n(t/z)$}
that are computed according to the recurrence relation
\begin{equation}
\hat{P}_n(t;z) = (2 - 1/n)t\hat{P}_{n-1}(t;z) - (1 - 1/n)z^2\hat{P}_{n-2}(t;z) \;,
\end{equation}
for
\mbox{$n = 1,2,\ldots$},
subject to the boundary values
\begin{equation}
\hat{P}_{-1}(t;z) = 1/z \;, \quad \hat{P}_0(t;z) = 1 \;.
\end{equation}
The scaled marginal integral
\begin{equation}
M_{\text{scl}}(n) \equiv \pi^{-1}(\alpha - \beta)^n M_n(\alpha,\beta) \;,
\end{equation}
in terms of which the cut contribution to the low-priority  marginal PMF is expressed as
\begin{equation}
P_{\text{cut}}(n) = \frac{(1-r)r^{n-1}}{2}{\cdot}M_{\text{scl}}(n) \;,
\label{pmargcut}
\end{equation}
is computed from the scaled Legendre polynomials as follows.
With
\begin{equation}
D \equiv z^2 -2\chi z + 1 = (\beta^2 - 1)/(\alpha^2 - 1) \;,
\end{equation}
we can write
\begin{equation}
M_{\text{scl}}(n) = -z\hat{P}_{n-1}(t;z) + (2\chi - z)\hat{P}_n(t;z)
     - \frac{\sqrt{D}}{z}\biggl[1 - \sqrt{D}\sum_{k=0}^{n}\hat{P}_k(t;z)\biggr] \;.
\label{Msc1}
\end{equation}
It is instructive to note the asymptotic behaviour
\begin{equation}
1 - \sqrt{D}\sum_{k=0}^n\hat{P}_k(t;z) \asym{n\to\infty}
     \sqrt{\frac{\beta+1}{\beta-1}}z\hat{P}_n(t;z) \;,
\end{equation}
provided
\mbox{$\beta \neq 1$}.
When
\mbox{$\beta = 1$},
the foregoing expression is identically unity for all
\mbox{$n < \infty$}
since
\mbox{$\sqrt{D} = 0$}.
This enables us to establish the bounding property
\begin{equation}
0 < \frac{\sqrt{D}}{z}\biggl[1 - \sqrt{D}\sum_{k=0}^n\hat{P}_k(t;z)\biggr] <
     \sqrt{\frac{\beta+1}{\beta-1}}{\cdot}\sqrt{D}\hat{P}_n(t;z) \;.
\label{bound}
\end{equation}

The generating function for the Legendre polynomials implies that
\begin{equation}
\sqrt{D}\sum_{k=0}^\infty \hat{P}_k(t;z) = 1\;.
\label{SumP}
\end{equation}
Use of this identity leads directly to an
alternative expression for the scaled marginal integral, given by
\begin{equation}
M_{\text{scl}}(n) = -z\hat{P}_{n-1}(t;z) + (2\chi - z)\hat{P}_n(t;z)
     - \frac{D}{z}\sum_{k=n+1}^\infty \hat{P}_k(t;z) \;.
\label{Mscl2}
\end{equation}
In the numerical implementation, we construct the last term of (\ref{Msc1})
and test whether it satisfies the bound given by (\ref{bound}).
If the test is passed, then the scaled marginal integral is evaluated via (\ref{Msc1}).
If it fails (due to arithmetic underflow),
then the scaled marginal integral is evaluated instead via (\ref{Mscl2}),
where the sum is computed by sequentially adding terms until the desired level of
convergence is attained.

\subsubsection{Joint Distribution}
\label{JointDistLeg}
To deal with the joint distribution,
we introduce scaled $Q$-functions according to
\mbox{$\hat{Q}^\ell_m(t;z) \equiv z^m Q^\ell_m(t/z)$}.
Let
\begin{equation}
A_m \equiv \frac{2m+1}{m+1}{\cdot}t \;, \quad
B_m^\ell \equiv \frac{m^2 - \ell^2}{m(m+1)}{\cdot}z^2 \;.
\end{equation}
Then,
\begin{equation}
\hat{Q}_{m+1}^\ell(t,z) = A_m \hat{Q}_{m}^\ell(t,z) - B_m^\ell \hat{Q}_{m-1}^\ell(t,z) \;,
\end{equation}
for
\mbox{$\ell = 0,1,\ldots$},
\mbox{$m = 1,2,\ldots$},
subject to the boundary values
\begin{equation}
\hat{Q}_0^\ell(t,z) = 1 \;, \quad \hat{Q}_1^\ell(t,z) = t + \ell z \;.
\end{equation}
We use the scaled $Q$-functions to construct the scaled joint integrals
\begin{equation}
J_{\text{scl}}(\ell,m) \equiv \pi^{-1}(-1)^\ell(\alpha - \beta)^{m+1}
     J_m^\ell(\alpha, \beta) \;,
\label{Jscl}
\end{equation}
as
\begin{equation}
\begin{split}
J_{\text{scl}}(\ell,m) &= d_{\ell+1}(\beta) +
     d_{\ell+1}(\alpha)\biggl[z\hat{Q}_{m}^{\ell+1}(\chi z,z) \\
     &\quad{}- (1+\chi)\biggl(
\sum_{k=0}^m \hat{Q}_k^\ell(\chi z, z) - \frac{1 - z/\chi}{1 + 1/\chi}
     \sum_{k=0}^m \hat{Q}_k^{\ell+1}(\chi z, z)
\biggr)\biggr] \;,
\end{split}
\end{equation}
and we may note the generalization of (\ref{SumP}), given by
\begin{equation}
\sqrt{D}\sum_{k=0}^\infty \hat{Q}_k^\ell(t;z) = \frac{d_\ell(\beta)}{d_\ell(\alpha)} \;.
\end{equation}
One can also confirm that
\mbox{$J_{\text{scl}}(\ell,0) =  d_{\ell+1}(\beta) -  d_{\ell+1}(\alpha)$}.

Given the definition of the scaled joint integrals in (\ref{Jscl}),
the cut contribution to the wait-conditional joint PMF becomes
\begin{equation}
P_{\text{cut}}(\ell,m) = (1-r)\sqrt{r_{\text{hi}}}\left[
     r_{\text{hi}}^{\ell/2}r^{m-1}
     J_{\text{scl}}(\ell,m-1)
     - r_{\text{hi}}^{(\ell-1)/2}r^m
     J_{\text{scl}}(\ell-1,m)\right] \;.
\end{equation}
As a sanity check, one can confirm, with the aid of (\ref{d1alpha}) and (\ref{d1beta}),
that this expression leads to the results
\begin{align}
\begin{aligned}
P_{\text{cut}}(0,0)    &= (1-r){\cdot}\min\left(1, r_{\text{hi}}/r^2\right) \,, \\
P_{\text{cut}}(0,m+1)  &= r_{\text{lo}}{\cdot}P_{\text{cut}}(m) \;, \\
P_{\text{cut}}(\ell,0) &= P_{\text{xhi}}(\ell) - P_{\text{pol}}(\ell,0) \;,
\end{aligned}
\end{align}
where the cut contribution to the low-priority marginal PMF
$P_{\text{cut}}(m)$ is given by (\ref{pmargcut}),
$P_{\text{pol}}(\ell,0)$ follows from (\ref{pjointpol})
and
\mbox{$P_{\text{xhi}}(\ell) \equiv P(\ell,0) = (1-r)\zeta_-^\ell(0)$} \citep{NP:Zuk23A}.
An intermediate step in deriving the last relationship is
\begin{equation}
P_{\text{cut}}(\ell,0) = (1-r)r_{\text{hi}}^{\ell/2}\left[d_\ell(\alpha) - d_\ell(\beta)
     + (\sqrt{r_{\text{hi}}}/r){\cdot}d_{\ell+1}(\beta)\right] \;,
\end{equation}
which is established with the aid of the identity
\begin{equation}
I_0^\ell(\gamma) + \gamma I_0^{\ell+1}(\gamma) = \pi(-1)^\ell d_{\ell+1}(\gamma) \;,
\end{equation}
with
\mbox{$\gamma = \alpha,\beta$}.
In fact, it may be noted that
\begin{equation}
P_{\text{pol}}(\ell,m) = (1-r)r^{m-1}r_{\text{hi}}^{\ell/2}\left[
     r d_\ell(\beta) - \sqrt{r_{\text{hi}}}d_{\ell+1}(\beta)\right] \;.
\end{equation}
Finally, we mention that the large-$\ell$ asymptotics of the scaled joint integrals are given by
\begin{equation}
J_{\text{scl}}(\ell,m) \asym{\ell\to\infty} d_{\ell+1}(\beta)
      - \frac{(\ell+1)^m}{m!}\left(\frac{\alpha-\beta}{\sqrt{\alpha^2-1}}\right)^m
      d_{\ell+1}(\alpha) \;.
\end{equation}
A derivation is presented in Appendix~\ref{JAsym}.
The large-$\ell$ behaviour of the joint PMF $P(\ell,m)$
for fixed $m$ follows directly from this.

\section{Numerical Results}
\label{Numerical}
In the present work, we have introduced two approaches to the practical
calculation of both the marginal and joint queue-length distributions.
We shall test the numerical performance of these two methods against each other
and against two previously developed methods,
namely, the quadratic-recurrence method from \citep{NP:Zuk23}
and the Fast Fourier Transform (FFT) method from \citep{NP:Zuk23A}.
We subject the algorithms to a number of tests,
each with an associated measure of performance (MOP) as follows:

\begin{figure}
\FIGURE
{\includegraphics[width=\wscl\linewidth, height=\hscl\linewidth]{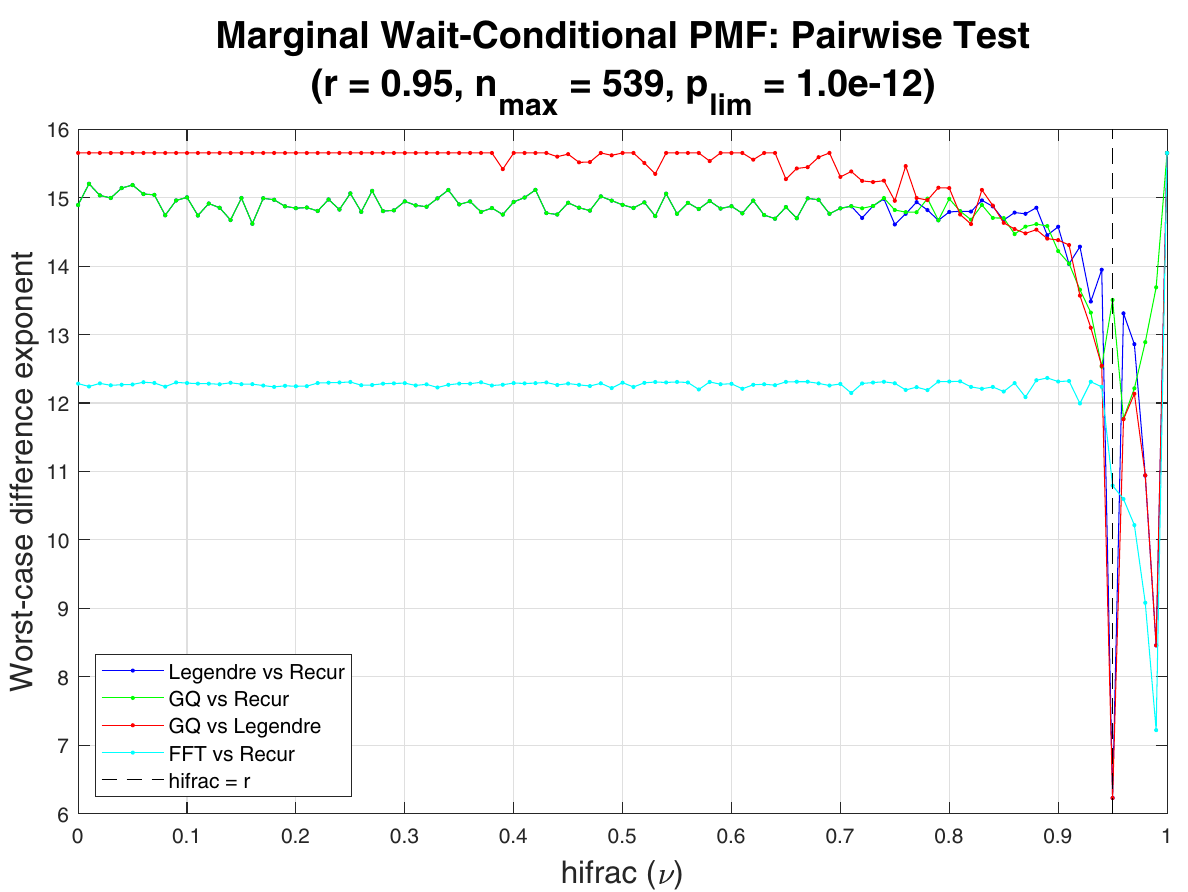}}
{\hphantom{x}\label{LoMargTest}}
{Pairwise comparisons among four algorithms for the low-priority marginal distribution
     of the queue length as a function of the
     fraction of high-priority arrivals $\nu$, with
     total traffic intensity $r = 0.95$.
     The values on the vertical axis indicate the number of decimal places of agreement.}
\end{figure}

\begin{figure}
\FIGURE
{\includegraphics[width=\wscl\linewidth, height=\hscl\linewidth]{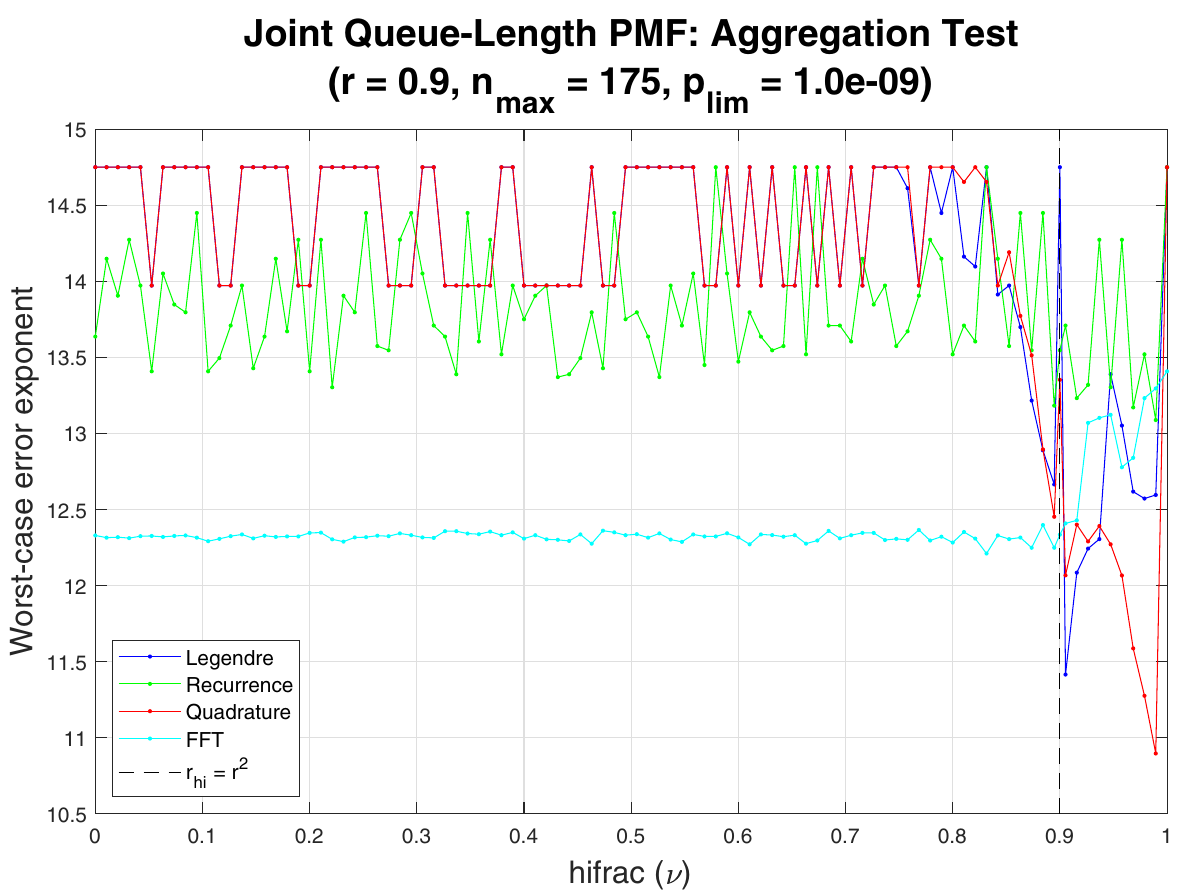}}
{\hphantom{x}\label{AggTest}}
{Aggregation tests for the joint probability distribution of the queue lengths as a function of the
     fraction of high-priority arrivals $\nu$, with
     total traffic intensity $r = 0.9$, comparing four algorithms as indicated.
     The values on the vertical axis indicate the number of decimal places of agreement.}
\end{figure}

\begin{figure}
\FIGURE
{\includegraphics[width=\wscl\linewidth, height=\hscl\linewidth]{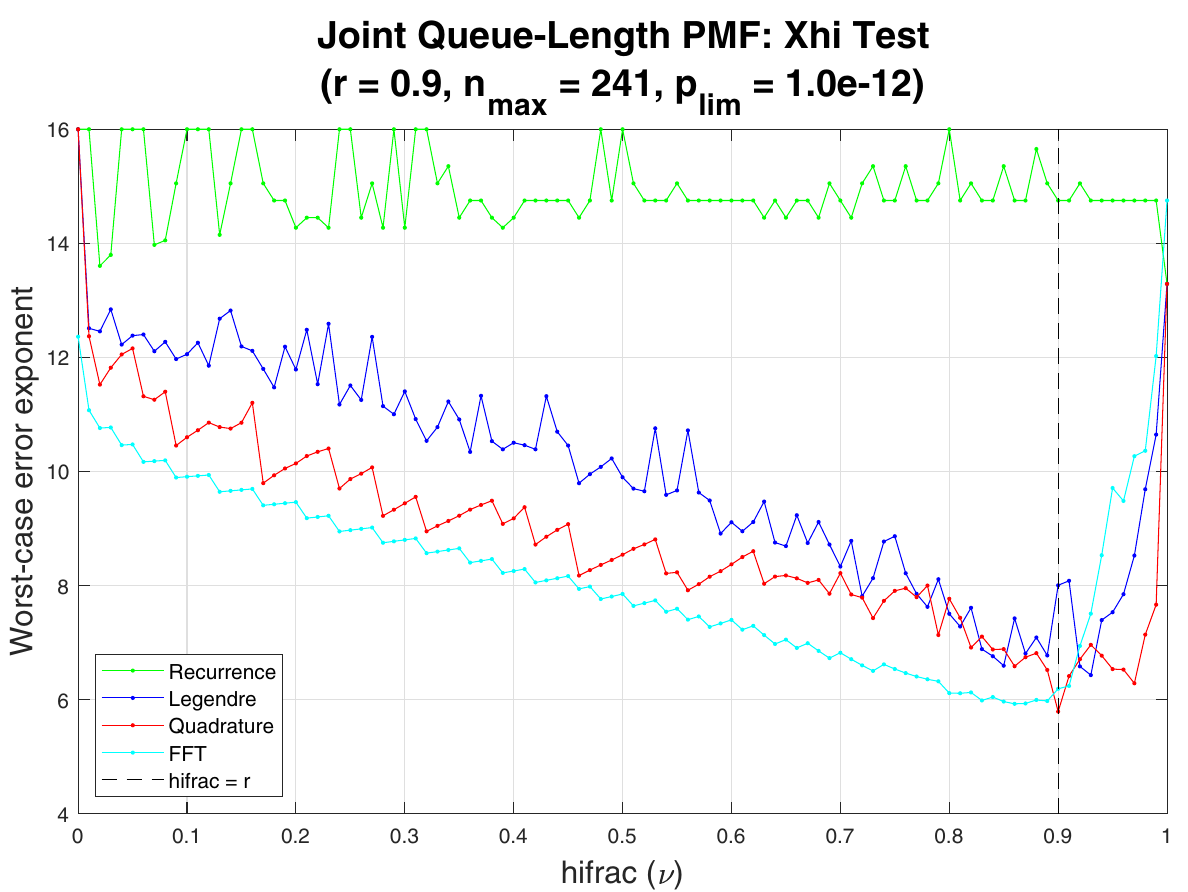}}
{\hphantom{x}\label{XhiTest}}
{Xhi-tests for the joint probability distribution of the queue lengths as a function of the
     fraction of high-priority arrivals $\nu$, with
     total traffic intensity $r = 0.9$, comparing four algorithms as indicated.
     The values on the vertical axis indicate the number of decimal places of agreement.}
\end{figure}

\begin{figure}
\FIGURE
{\includegraphics[width=\wscl\linewidth, height=\hscl\linewidth]{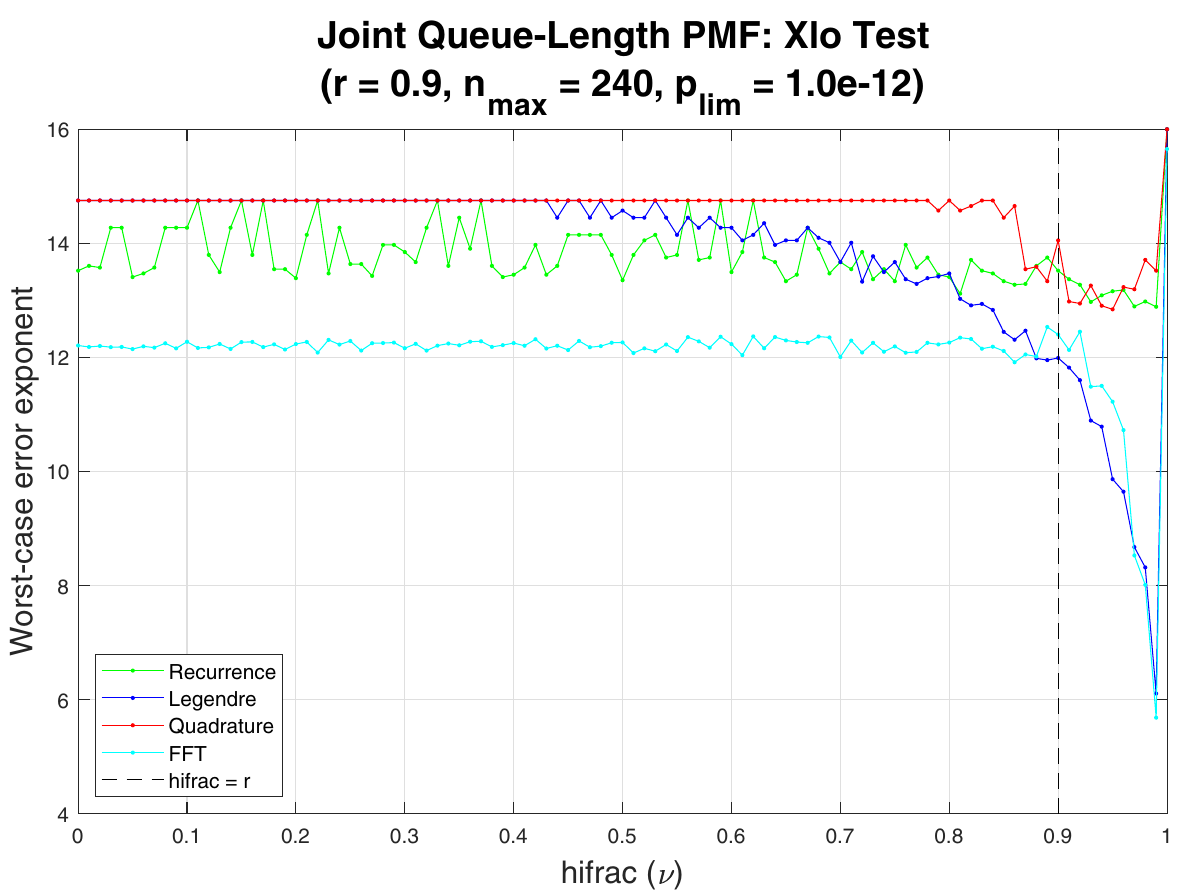}}
{\hphantom{x}\label{XloTest}}
{Xlo-tests for the joint probability distribution of the queue lengths as a function of the
     fraction of high-priority arrivals $\nu$, with
     total traffic intensity $r = 0.9$, comparing four algorithms as indicated.
     The values on the vertical axis indicate the number of decimal places of agreement.}
\end{figure}

\paragraph{Pairwise Test:}
We compare the low-priority marginal PMFs computed via the various methods with
each other in a pairwise manner.
To compare method A with method B, we adopt the MOP
\begin{equation}
\Xi_{\text{pair}} \equiv -\max_{m \geq 0}\left\{\log_{10}\left(|\ln(P_{\text{lo}}^{(\text{A})}(m)) -
    \ln{(}P_{\text{lo}}^{(\text{B})}(m))|\right)\right\} \;,
\end{equation}
where the maximum is taken over all values
\mbox{$0 \leq m \leq n_{\text{lim}}$}
such that
\begin{equation}
\min\{P_{\text{lo}}^{(\text{A})}(m), P_{\text{lo}}^{(\text{B})}(m)\} > p_{\text{lim}} > 0 \;.
\end{equation}
Since we are working in double-precision arithmetic\footnote{All computation is performed in {\sc Matlab} R2020a,
which implements IEEE Standard 754 for double precision.},
all MOPs of this kind are capped
at a maximum allowed value of $16$.
The interpretation of $\Xi_{\text{pair}}$
(and similarly for all of the subsequent MOPs)
is that it indicates the number of decimal places of numerical agreement
in the worst case.
Finite values for the limiting probability level $p_{\text{lim}}$, typically around  $10^{-12}$,
are adopted as one cannot expect any digital algorithm to work infinitely far into the tail.
Finite values of the limiting queue occupancy $n_{\text{lim}}$ are adopt to keep computational
load under control.

\paragraph{Aggregation Test:}
The aggregated queue-length distribution describes the total number of entities in the queue,
regardless of priority level.
The exact result for the aggregated PMF is \citep{NP:Zuk23,NP:Zuk23A}
\begin{equation}
P_{\text{agg}}^{\text{(ex)}}(k) = (1-r)r^k \;,
\label{Pagg}
\end{equation}
for
\mbox{$k = 0,1,\ldots$}.
We construct the aggregated probabilities from the numerically computed
joint PMF according to
\begin{equation}
P_{\text{agg}}(k) = \sum_{m=0}^k P(k-m,m)
\end{equation}
and test against the exact result using the MOP
\begin{equation}
\Xi_{\text{agg}} \equiv -\max_{k \geq 0}\left\{\log_{10}\left(|\ln(P_{\text{agg}}(k)) -
    \ln{(}P_{\text{agg}}^{(\text{ex})}(k))|\right)\right\} \;,
\end{equation}
where the maximum is taken over all values
\mbox{$0 \leq k \leq n_{\text{lim}}$}
such that
\mbox{$P_{\text{agg}}^{(\text{ex})}(k) > p_{\text{lim}} > 0$}.

\paragraph{Xhi-Test:}
The exact result for the exclusively-high queue-length probabilities is \cite{NP:Zuk23A}
\begin{equation}
P_{\text{xhi}}^{\text{(ex)}}(\ell) \equiv P(\ell,0) = (1-r){\cdot}
     \left[\frac{1 + r - \sqrt{(1+r)^2 - 4r_{\text{hi}}}}{2}\right]^\ell \;,
\label{Pxhi}
\end{equation}
for
\mbox{$\ell = 0,1,\ldots$}.
We construct the exclusively-high probabilities from the numerically computed
joint PMF according to
\mbox{$P_{\text{xhi}}(\ell) = P(\ell,0)$}
and test against the exact result using the MOP
\begin{equation}
\Xi_{\text{xhi}} \equiv -\max_{\ell \geq 0}\left\{\log_{10}\left(|\ln(P_{\text{xhi}}(\ell)) -
    \ln{(}P_{\text{xhi}}^{(\text{ex})}(\ell))|\right)\right\} \;,
\end{equation}
where the maximum is taken over all values
\mbox{$0 \leq \ell \leq n_{\text{lim}}$}
such that
\mbox{$P_{\text{xhi}}^{(\text{ex})}(\ell) > p_{\text{lim}} > 0$}.

\paragraph{Xlo-Test:}
An exact relationship between that exclusively-low queue probabilities and
the low-priority marginal PMF is given by \citep{NP:Zuk23,NP:Zuk23A}
\begin{equation}
P_{\text{xlo}}^{(\text{ex})}(m) \equiv P(0,m)
     = (1-r)\delta_{m0} + (1 - \delta_{m0})r_{\text{lo}}{\cdot}P_{\text{lo}}(m-1) \;,
\label{xlotest}
\end{equation}
for
\mbox{$m = 0,1,2,\ldots$},
where we can formally set
\mbox{$P_{\text{lo}}(-1) \equiv 0$}.
We construct the exclusively-low probabilities from the numerically computed
joint PMF according to
\mbox{$P_{\text{xlo}}(m) = P(0,m)$}
and test against the low-priority marginal using the MOP
\begin{equation}
\Xi_{\text{xlo}} \equiv -\max_{m > 0}\left\{\log_{10}\left(|\ln(P_{\text{xlo}}(m)) -
    \ln{(}r_{\text{lo}}P_{\text{lo}}(m-1))|\right)\right\} \;,
\end{equation}
where the maximum is taken over all values
\mbox{$0 < m \leq n_{\text{lim}}$}
such that
\mbox{$P_{\text{xlo}}^{(\text{ex})}(m) > p_{\text{lim}} > 0$}.

In Figures~\ref{LoMargTest}--\ref{XloTest}, the MOP values relevant to each test
are displayed on the vertical axis against the full range of
high-priority arrival fraction (hifrac)
\mbox{$0 \leq \nu \leq 1$}
on the horizontal axis.
For the marginal distribution,
Figure~\ref{LoMargTest} presents results of the pairwise tests for
\mbox{$r = 0.95$}
with
\mbox{$p_{\text{lim}} = 10^{-12}$}.
The blue curve compares the Legendre polynomial method of Section~\ref{MargDistLeg}
and the quadratic recurrence method from \citep{NP:Zuk23}.
The green curve compares the Gaussian quadrature method obtained from (\ref{MargPol}) and (\ref{PGQ}) of
Section~\ref{Marginal} and the aforementioned recurrence method.
The red curve compares the quadrature and Legendre methods.
The cyan curve compares the recurrence method and the FFT method from \citep{NP:Zuk23A}.

Here, and subsequently, $n_{\text{max}}$ denotes the maximum queue-length
actually sampled in achieving the limiting probability level $p_{\text{lim}}$.
The FFT method is configured to a relative error of $10^{-12}$,
at least for small values of $\nu$,
and this is reflected in all of the graphs.
The fact that the blue and red curves dip in accuracy at the critical point
(\mbox{$\nu = r$}),
but the green curve does not, indicates that only the Legendre method struggles
in this region, while the quadratic-recurrence and quadrature methods
are able to maintain accuracy.
The fact that the blue and green curves (where the quadratic recurrence is the common method)
align closely but lie below the red curve suggests that the quadrature and Legendre methods are
slightly more accurate, at least away from the critical point
(indicated by the dashed black line).

In Figures~\ref{AggTest}--\ref{XloTest}, that test the joint distribution,
the red curves correspond to the Gaussian quadrature method of Section~\ref{Joint}.
The blue curves correspond to the Legendre polynomial method of Section~\ref{JointDistLeg}.
The green curves result from the quadratic recurrence method from \citep{NP:Zuk23}, and
the cyan curves from the FFT method described in \citep{NP:Zuk23A}.
Figure~\ref{AggTest} presents results of the aggregation test for
\mbox{$r = 0.9$}
with
\mbox{$p_{\text{lim}} = 10^{-9}$}.
Figure~\ref{XhiTest} presents results of the xhi-test for
\mbox{$r = 0.9$}
with
\mbox{$p_{\text{lim}} = 10^{-12}$}.
Figure~\ref{XloTest} presents results of the xlo-test for
\mbox{$r = 0.9$}
with
\mbox{$p_{\text{lim}} = 10^{-12}$}.
All methods exhibit good levels of numerical performance.

\begin{figure}
\FIGURE
{\includegraphics[width=\wscl\linewidth, height=\hscl\linewidth]{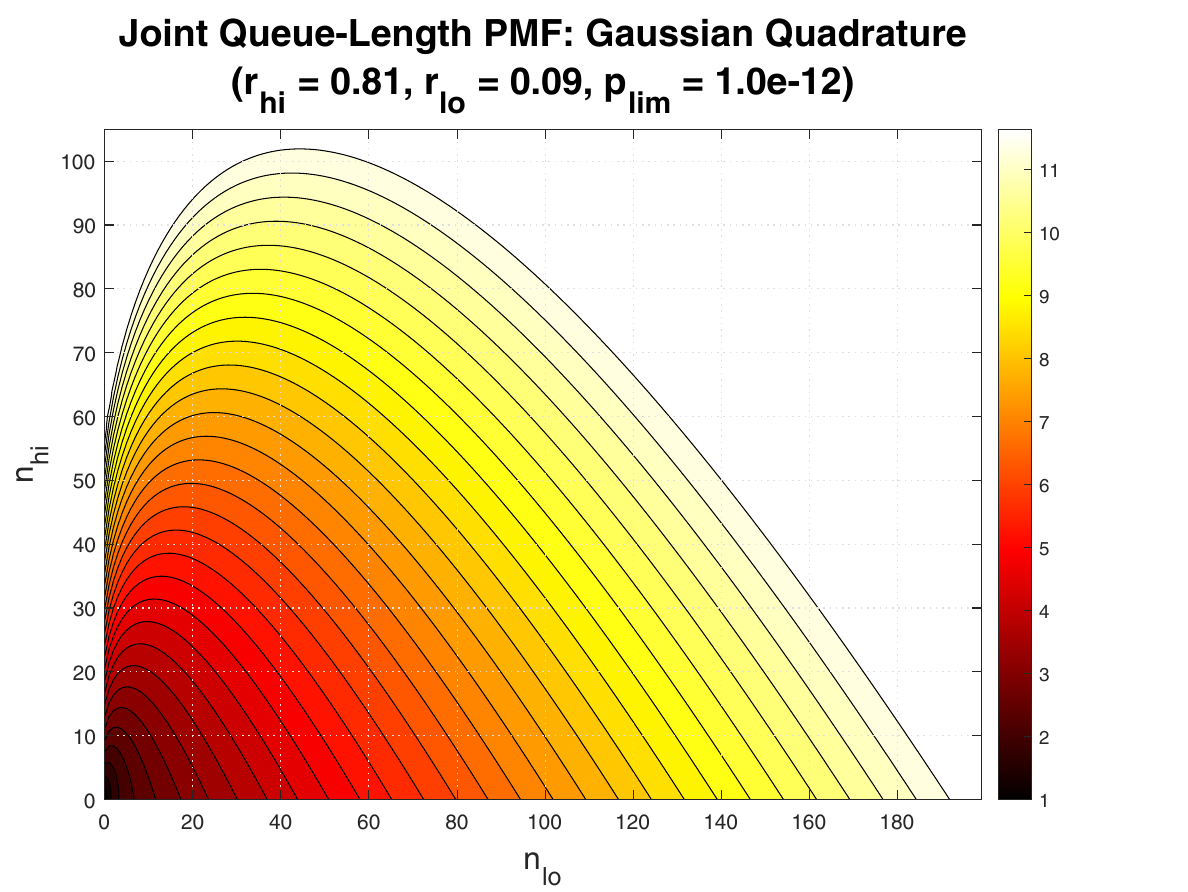}}
{\hphantom{x}\label{GQTest1}}
{Two-dimensional level-set mapping for the joint probability distribution of the queue lengths
     for total traffic intensity $r = 0.9$ and fraction of high-priority arrivals
     at the critical point $\nu = r$.}
\end{figure}

\begin{figure}
\FIGURE
{\includegraphics[width=\wscl\linewidth, height=\hscl\linewidth]{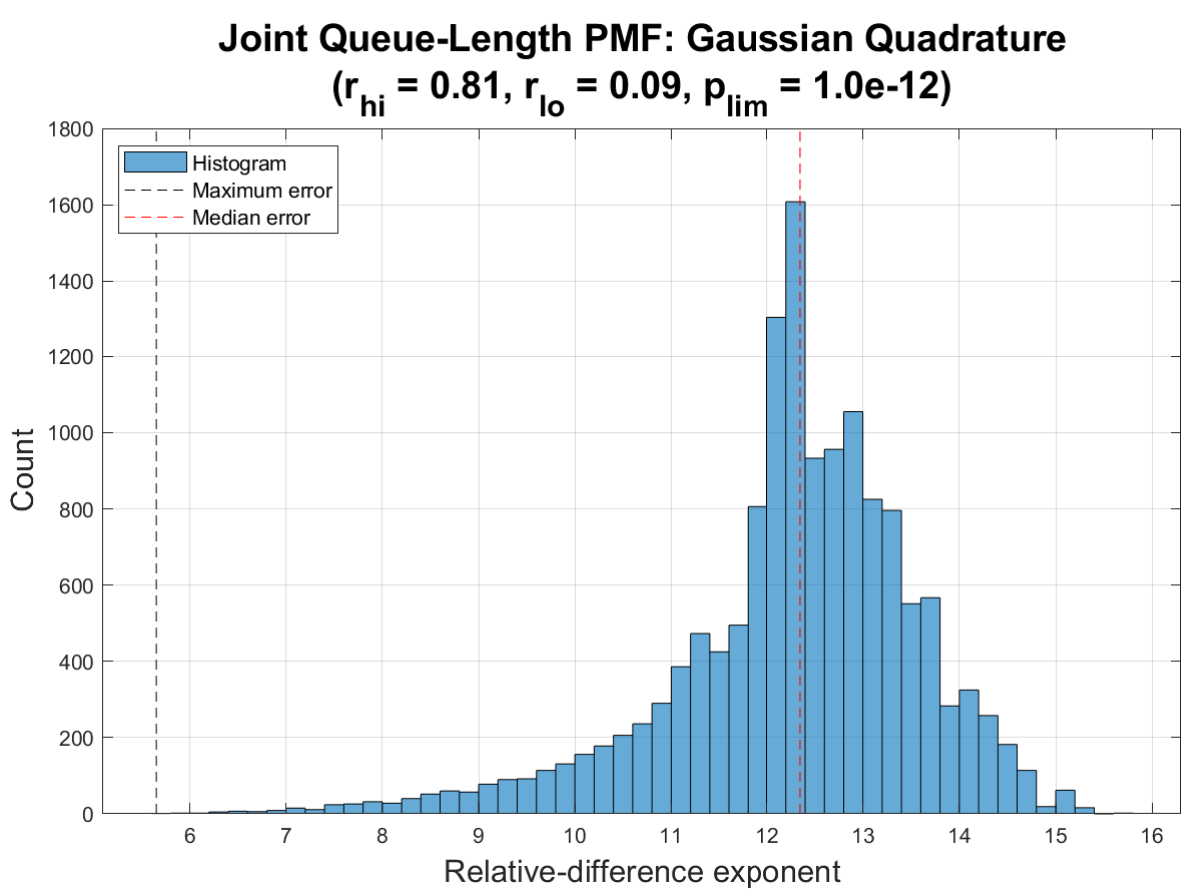}}
{\hphantom{x}\label{GQTest2}}
{Error histogram, comparing the quadrature and quadratic-recurrence methods,
     for the joint probability distribution of the queue lengths
     for total traffic intensity $r = 0.9$ and fraction of high-priority arrivals
     at the critical point $\nu = r$.}
\end{figure}

\begin{figure}
\FIGURE
{\includegraphics[width=\wscl\linewidth, height=\hscl\linewidth]{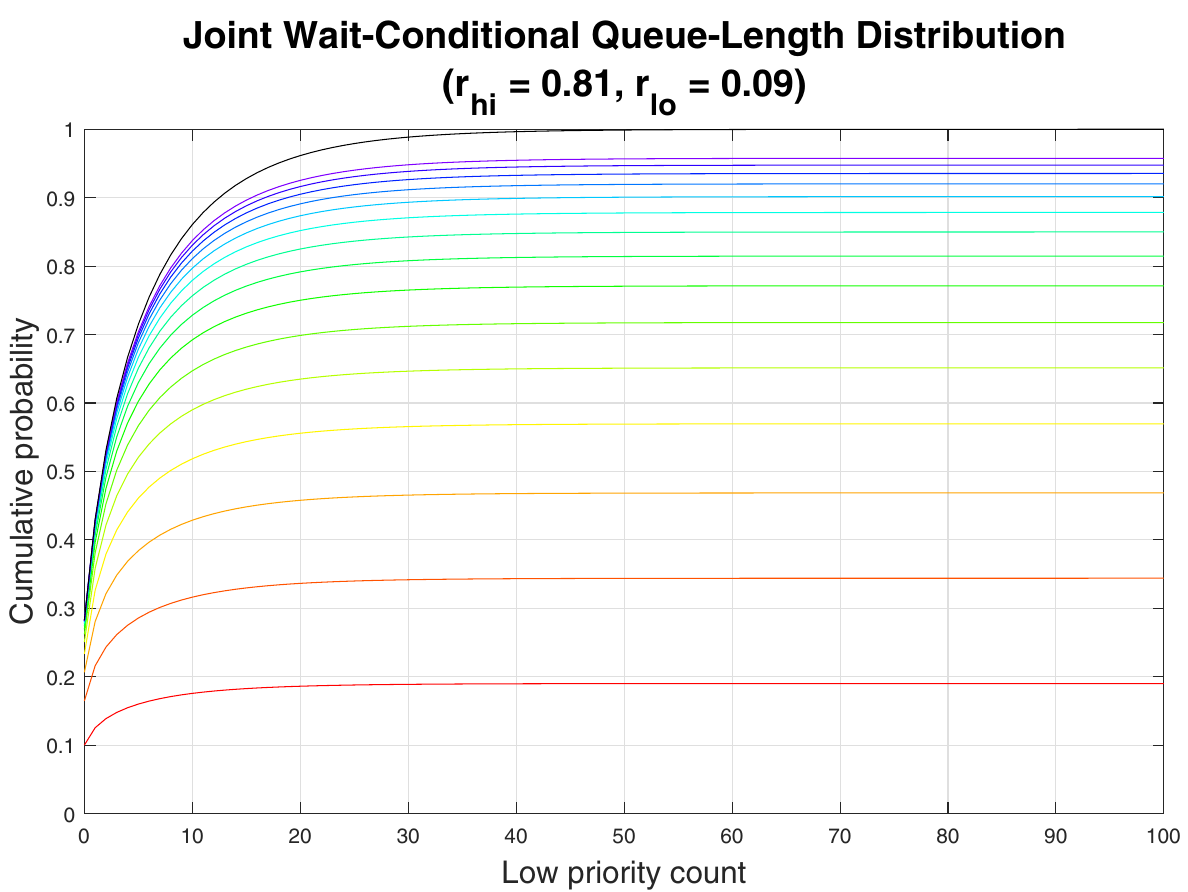}}
{\hphantom{x}\label{CDFSect}}
{CDF sections for the joint probability distribution of the queue lengths
     for total traffic intensity $r = 0.9$ and fraction of high-priority arrivals
     at the critical point $\nu = r$.}
\end{figure}

In order to illustrate the qualitative features of the queue-length distribution,
in Figure~\ref{GQTest1}, we present a two-dimensional contour plot with level curves
for the joint queue-length PMF computed via the quadrature method at the critical point
(\mbox{$\nu = r$})
for
\mbox{$r = 0.9$}.
This is computationally the worst-case scenario for a given total traffic intensity $r$.
The level curves are displayed on a (base-$10$) logarithmic scale down to a limiting
probability level of
\mbox{$p_{\text{lim}} = 10^{-12}$},
with the low-priority occupancy on the horizontal axis and the high-priority
occupancy on the vertical axis.
We also benchmark this case against the quadratic-recurrence method, with the results
for the number of decimal places of agreement
presented as a histogram in Figure~\ref{GQTest2}.
Finally, in Figure~\ref{CDFSect},
we plot one-dimensional sections of the cumulative distribution function (CDF)
for various constant values of the high-priority queue occupancy, {\it i.e.}
\begin{equation}
F_\ell(n) \equiv \sum_{m=0}^n P(\ell,m) \;.
\end{equation}
The curves correspond to high-priority queue-lengths
\mbox{$\ell = 0,1,\ldots\,14$},
with the uppermost black curve representing
\mbox{$\ell = 100$}.

\section{Conclusions}
\label{Concl}
The analytic approach to the calculation of the joint and marginal queue-length distributions
for the two-level non-preemptive priority queue studied here
leads to integral representations that are
amenable to both efficient numerical quadrature rules and exact expressions as finite sums
of easily computable special functions that generalize the Legendre polynomials.
One benefit of integral forms obtained here
is that the large queue-length asymptotic behaviour becomes apparent.
Another is that there is scope for extending the methodology to cater for additional complications,
such as multiple arrival classes.

\begin{APPENDICES}
\section{No-Wait Probability}
\label{NoWait}
Let $c$ denote the number of servers, and let $\mathscr{N}_{\text{sys}}$ be the
RV for system occupancy. It is well-known that
\mbox{$\Pr(\mathscr{N}_{\text{sys}} = n) = Ar^n$}
for all
\mbox{$n \geq c$},
for some constant $A$
and, by definition, we have
\mbox{$P_{\text{NW}} = \Pr(\mathscr{N}_{\text{sys}} \leq c-1)$}
for the no-wait probability.
Therefore,
\begin{align}
\begin{aligned}
P_{\text{NW}} &= \sum_{n=0}^{c-1} \Pr(\mathscr{N}_{\text{sys}} = n) \\
     &= 1 - \sum_{n=c}^\infty \Pr(\mathscr{N}_{\text{sys}} = n) \\
     &= 1 - Ar^c/(1-r) \\
     &= 1 - \Pr(\mathscr{N}_{\text{sys}} = c)/(1-r) \;.
\end{aligned}
\end{align}
Furthermore,
\mbox{$\Pr(\mathscr{N}_{\text{sys}} = c) =  [(rc)^c/c!]\Pr(\mathscr{N}_{\text{sys}} = 0)$},
and \citep{NP:Gnedenko89}
\begin{equation}
\frac{1}{\Pr(\mathscr{N}_{\text{sys}} = 0)} = \sum_{n=0}^c \frac{(rc)^n}{n!}
     + \frac{(rc)^c}{c!(1/r-1)} \;.
\end{equation}
The result for $P_{\text{NW}}$ follows, and may be compactly stated as
\begin{equation}
\frac{1}{1-P_{\text{NW}}} = 1 + \frac{(1-r)c!}{(rc)^c}
     \sum_{n=0}^{c-1}\frac{(rc)^n}{n!} \;.
\end{equation}

\section{Laguerre Series for the Waiting Time}
\label{Laguerre}
Let us write
\mbox{$p_n \equiv P_\mathscr{L}(n)$}
so that
\begin{equation}
p_n = \frac{1}{n!}{\cdot}\left.\frac{d^n}{dp^n} g_{\text{lo}}(z)\right|_{z=0}
     = \oint_{\mathcal{C}_\eta}\frac{dz}{2\pi i}\, \frac{g_{\text{lo}}(z)}{z^{n+1}} \;,
\end{equation}
where the closed anti-clockwise traversed integration contour $\mathcal{C}_\eta$, initially taken to be a
circle centred on the origin of radius
\mbox{$0 < \eta < 1/r$},
is deformed into  the Bromwich contour $\Gamma$ on which
\mbox{$z = \delta +iy$}
with
\mbox{$-\infty < y < +\infty$},
for some
\mbox{$1 < \delta < 1/r$}
so that
\mbox{$\Rez z > 1$}.
Then, we consider the Binomial transform of the $p_n$ sequence
\begin{align}
\begin{aligned}
\gamma_n &\equiv \sum_{k=0}^n \binom{n}{k}(-1)^k p_k \\
     &= \int_\Gamma\frac{dz}{2\pi i}\, g_{\text{lo}}(z)\sum_{k=0}^n \binom{n}{k}(-1)^k
     \frac{1}{z^{k+1}} \\
     &= \int_\Gamma\frac{dz}{2\pi i}\, \frac{g_{\text{lo}}(z)}{z}{\cdot}
     (1-1/z)^n \;.
\end{aligned}
\label{gamman}
\end{align}
We now recall the generating function for the Laguerre polynomials $L_n(x)$, given by
\begin{equation}
\sum_{n=0}^\infty \tau^n L_n(x) = \frac{1}{1-\tau}e^{-\tau x/(1-\tau)} \;.
\end{equation}
On setting
\mbox{$\tau = 1 - 1/z$}
and
\mbox{$x = r_{\text{lo}}t$},
this reads
\begin{equation}
\sum_{n=0}^\infty \left(1 - 1/z\right)^nL_n(r_{\text{lo}}t) = z e^{-r_{\text{lo}}t(z-1)} \;,
\end{equation}
from which it follows that
\begin{equation}
\sum_{n=0}^\infty \gamma_n L_n(r_{\text{lo}}t) = \int_\Gamma\frac{dz}{2\pi i}\,
     g_{\text{lo}}(z) e^{-r_{\text{lo}}t(z-1)} \;.
\label{Lsum}
\end{equation}
Now, (\ref{DistLittle}) can be written as
\begin{equation}
P_\mathscr{L}(k) = \int_0^\infty dt'\, P'_\mathscr{W}(t') \frac{t^{\prime k}}{k!}e^{-t'} \;,
\end{equation}
where
\mbox{$t' = r_{\text{lo}}t$}
and
\mbox{$P'_\mathscr{W}(t')dt' = P_\mathscr{W}(t)dt$}.
We also have
\begin{equation}
L_n(t) = \sum_{k=0}^n \binom{n}{k}\frac{(-1)^k}{k!} t^k \;.
\end{equation}
These may be combined to yield
\begin{equation}
\int_0^\infty dt'\, e^{-t'}L_n(t') P'_\mathscr{W}(t') =
     \sum_{k=0}^n \binom{n}{k}(-1)^k \int_0^\infty dt'\,  P'_\mathscr{W}(t')
     \frac{t^{\prime k}}{k!}e^{-t'} = \gamma_n \;.
\end{equation}
The completeness relation \citep{NP:Li08}
\begin{equation}
\sum_{n=0}^\infty L_n(t)L_n(t') = e^{(t+t')/2}{\cdot}\delta(t-t') \;,
\end{equation}
\mbox{$\delta(t)$}
being the Dirac delta function,
implies the convergent Laguerre-series representation
\begin{equation}
P'_\mathscr{W}(t') = \sum_{n=0}^\infty \gamma_n L_n(t') \;,
\label{LSeries}
\end{equation}
and we conclude, via (\ref{Lsum}),  that
\begin{equation}
P_\mathscr{W}(t)dt = r_{\text{lo}}dt{\cdot}\int_\Gamma\frac{dz}{2\pi i}\,
     g_{\text{lo}}(z) e^{-r_{\text{lo}}t(z-1)} \;.
\end{equation}
The Bromwich contour $\Gamma$, with
\mbox{$\Re z > 1$},
may now be deformed as
\mbox{$\Gamma\mapsto \mathcal{C}_{\text{pol}} \cup \mathcal{C}_{\text{cut}}$}
so that it wraps around the pole and cut of $g_{\text{lo}}(z)$.
This serves to reproduce (\ref{WaitContour}).

The convergence claim for (\ref{LSeries}) follows from the theorem that
the formal expansion
\begin{equation}
\gamma(t) \equiv \sum_{k=0}^\infty \gamma_k L_k(t)
\end{equation}
converges in $L_2$ if and only if
\mbox{$\|\gamma\|^2_{L_2} < \infty$},
where
\begin{equation}
\|\gamma\|^2_{L_2} \equiv \int_0^\infty dt\, e^{-t}|\gamma(t)|^2
     = \sum_{k=0}^\infty |\gamma_k|^2 \;.
\end{equation}
In the present application, we have
\mbox{$\gamma(t) = P'_{\mathscr{W}}(t)$}
for which, trivially,
\begin{equation}
\|P'_\mathscr{W}\|_{L_2}^2 = \int_0^\infty dt\, e^{-t}\left[P'_\mathscr{W}(t)\right]^2
     < \int_0^\infty dt\, e^{-t} = 1 \;.
\end{equation}

One should note that, while the coefficients $\gamma_n$ of the Laguerre series
(\ref{LSeries}) are constructed as (potentially)
numerically problematic alternating combinatorial sums, the last equality in
(\ref{gamman}) shows that they can be efficiently computed by means of the very same
iterative quadrature method derived for the queue-length probabilities $p_n$.
Explicitly,
\begin{equation}
\begin{split}
\gamma_n &= \Theta(r^2 - r_{\text{hi}})
     \left[1 - \frac{r(1-r)}{r_{\text{lo}}}\right]{\cdot}\frac{(1-r)^{n+1}}{r} \\
     &\quad {}+ \frac{1-r}{2\pi r}\int_0^1 du\,\sqrt{u(1-u)}{\cdot}
          \frac{(u+c)^n}{(u + a)^{n+1}(u + b)} \;,
\end{split}
\end{equation}
where
\mbox{$c \equiv (x_- - 1)/x_{\text{dif}} = (1-r)a + rb$}.

\section{Joint Asymptotics}
\label{JAsym}
To determine the large-$\ell$ behaviour of the joint integral (\ref{JIntg}),
we use the symmetry properties and periodicity of the integrand
to write it as an integral over the interval
\mbox{$\theta \in [0,2\pi)$}.
We then promote the integration
variable $\theta$ to the complex plane, and deform the contour from
\mbox{$\theta = 0$}
to
\mbox{$\theta = 2\pi$}
into three straight legs: (i) from $0$ to $0 + i\infty$,
(ii) from $0 + i\infty$ to $2\pi + i\infty$, and
(iii) from  $2\pi + i\infty$ to $2\pi$,
making small excursions on the vertical legs around the points
\mbox{$\theta = \acos(\alpha), \acos(\beta)$}.
The contribution from the horizontal contour vanishes,
while the vertical contours cancel apart from the excursions around the singular points.
As a result, we obtain the decomposition
\mbox{$J_{\text{scl}}(\ell,m) = J_{\text{scl}}^\alpha(\ell,m) + J_{\text{scl}}^\beta(\ell,m)$}
where, on writing
\mbox{$\theta = it$},
we have
\begin{equation}
J_{\text{scl}}^{\alpha,\beta}(\ell,m) = (\beta-\alpha)^{m+1}\oint_{C_{\alpha,\beta}}
     \frac{dt}{2\pi i}\, \frac{\sinh(t)e^{-\ell t}}{(\cosh(t)-\alpha)^{m+1}(\cosh(t) - \beta)} \;,
\end{equation}
where $C_\alpha$ is an infinitesimal anti-clockwise circle about the point
\mbox{$t = \acosh(\alpha)$},
and likewise for $C_\beta$ about the point
\mbox{$t = \acosh(\beta)$}.

It is clear that
\mbox{$J_{\text{scl}}^\beta(\ell,m) = d_{\ell+1}(\beta)$}.
To evaluate $J_{\text{scl}}^\alpha(\ell,m)$,
we set
\mbox{$\alpha \equiv \cosh(s)$}
so that
\begingroup
\addtolength{\jot}{0.5em}
\begin{align}
\begin{aligned}
J_{\text{scl}}^\alpha(\ell{-}1,m) &= \left(\frac{\beta-\alpha}{2}\right)^{m+1}\!\!e^{-\ell s}\oint_{C_\alpha}
     \frac{dt}{2\pi i}\, \frac{e^{-\ell (t-s)}}{\sinh^{m+1}((t-s)/2)} \\
     &\quad {}\times \frac{\sinh(t)}{\sinh^{m+1}((t+s)/2)[\cosh(t)-\beta]} \\
&= \left(\frac{\beta-\alpha}{2}\right)^{m+1}\frac{e^{-\ell s}}{\ell}\oint_{C_0}
      \frac{dt}{2\pi i}\,\frac{e^{-t}}{\sinh^{m+1}(t/(2\ell))} \\
     &\quad {}\times\frac{\sinh(s+t/\ell)}{\sinh^{m+1}(s+t/(2\ell))[\cosh(s+t/\ell)-\beta]} \;.
\end{aligned}
\end{align}
\endgroup
Taking the limit yields
\begin{equation}
J_{\text{scl}}^\alpha(\ell{-}1,m)
     \asym{\ell\to\infty} \frac{(\beta-\alpha)^{m+1}\ell^me^{-\ell s}}{\sinh^m(s)[\cosh(s)-\beta]}
     \oint_{C_0}\frac{dt}{2\pi i}\, \frac{e^{-t}}{t^{m+1}}
     = -\frac{\ell^m}{m!}\left(\frac{\alpha-\beta}{\sqrt{\alpha^2-1}}\right)^m d_\ell(\alpha) \,,
\end{equation}
where $\mathcal{C}_0$ is an infinitesimal anti-clockwise circle about the origin.
The result in the main text follows.
A by-product of this discussion is an alternative representation of the joint integral,
given by
\begin{equation}
J_{\text{scl}}(\ell,m) = d_{\ell+1}(\beta) + \oint_{C_{\alpha}}\frac{dw}{2\pi i}\,
     \frac{d_{\ell+1}(w)}{w-\beta}\left(\frac{\beta-\alpha}{w-\alpha}\right)^{m+1} \;.
\end{equation}
\end{APPENDICES}


%
\section*{Acknowledgments}
The authors gratefully acknowledge useful discussions with Dr.~Stephen Bocquet.


\bibliographystyle{informs2014} 
\bibliography{AnalyticNPQ} 

\begin{thebibliography}{32}
\providecommand{\natexlab}[1]{#1}
\providecommand{\url}[1]{\texttt{#1}}
\providecommand{\urlprefix}{URL }

\bibitem[{Almehdawe et~al.(2013)Almehdawe, Jewkes, \protect\BIBand{} {Q.-M.
  He}}]{NP:Almehdawe13}
Almehdawe E, Jewkes B, {Q-M He} (2013) A {M}arkovian queueing model for
  ambulance offload delays. \emph{European Journal of Operational Research}
  226(3):602--614.

\bibitem[{Bertsimas \protect\BIBand{} Mourtzinou(1997)}]{NP:Bertsimas97}
Bertsimas D, Mourtzinou G (1997) Multiclass queueing systems in heavy traffic:
  An asymptotic approach based on distributional and conservation laws.
  \emph{Operations Research} 45(3):470--487.

\bibitem[{Bertsimas \protect\BIBand{} Nakazato(1995)}]{NP:Bertsimas95}
Bertsimas D, Nakazato D (1995) The distributional {L}ittle's law and its
  applications. \emph{Operations Research} 43(2):298--310.

\bibitem[{Chawla(1970)}]{NP:Chawla70}
Chawla M (1970) Estimation of errors of {G}auss-{C}hebyshev quadratures.
  \emph{The Computer Journal} 13(1):107--109.

\bibitem[{Cohen(1956)}]{NP:Cohen56}
Cohen J (1956) Certain delay problems for a full availability trunk group
  loaded by two traffic sources. \emph{Philips Telecommunications Review}
  16(3):105--113.

\bibitem[{Davis(1966)}]{NP:Davis66}
Davis R (1966) Waiting time distribution of a multi-server, priority queuing
  sytem. \emph{Operations Research} 14(1):133--136.

\bibitem[{Erd\'{e}lyi(1940)}]{NP:Erdelyi40}
Erd\'{e}lyi A (1940) Some confluent hypergeometric functions of two variables.
  \emph{Proceedings of the Royal Society of {E}dinburgh} 60(3):344--361.

\bibitem[{Gail et~al.(1988)Gail, Hantler, \protect\BIBand{} Taylor}]{NP:Gail88}
Gail H, Hantler S, Taylor B (1988) Analysis of non-preemptive priority
  multi-server queue. \emph{Advances in Applied Probability} 20(4):852--879.

\bibitem[{Gautschi(2002)}]{NP:Gautschi02}
Gautschi W (2002) {G}auss quadrature approximations to hypergeometric
  functionsand confluent hypergeometric functions. \emph{Journal of
  Computational and Applied Mathematics} 139(1):173--187.

\bibitem[{Gautschi(2004)}]{NP:Gautschi04}
Gautschi W (2004) \emph{Orthogonal Polynomials: Computation and Approximation}
  (Oxford, UK: Oxford University Press), first edition.

\bibitem[{Gautschi(2005)}]{NP:Gautschi05}
Gautschi W (2005) Orthogonal polynomials (in {M}atlab). \emph{Journal of
  Computational and Applied Mathematics} 178(1):425--434.

\bibitem[{Gnedenko \protect\BIBand{} Kovalenko(1989)}]{NP:Gnedenko89}
Gnedenko B, Kovalenko I (1989) \emph{Introduction to Queueing Theory} (Boston,
  MA, USA: Birkh{\"{a}}user), second edition.

\bibitem[{Gradshteyn \protect\BIBand{} Ryzhik(2007)}]{NP:Gradshteyn07}
Gradshteyn I, Ryzhik I (2007) \emph{Table of Integrals, Series, and Products}
  (Burlington, MA, USA: Academic Press), seventh edition.

\bibitem[{Holley(1954)}]{NP:Holley54}
Holley J (1954) Waiting line subject to priorities. \emph{Operations Research}
  2(3):341--343.

\bibitem[{Hou \protect\BIBand{} Zhao(2020)}]{NP:Hou20}
Hou J, Zhao X (2020) Using a priority queueing approach to improve emergency
  department performance. \emph{Journal of Management Analytics} 7(1):28--43.

\bibitem[{Kao \protect\BIBand{} Narayanan(1990)}]{NP:Kao90}
Kao E, Narayanan K (1990) Computing steady state probabilities of a
  non-preemptive priority queue. \emph{{INFORMS} Journal on Computing}
  2(3):211--218.

\bibitem[{Karlin \protect\BIBand{} McGregor(1958)}]{NP:Karlin58}
Karlin S, McGregor J (1958) Many server queueing processes with {P}oisson input
  and exponential service times. \emph{Pacific Journal of Mathematics}
  8(1):87--118.

\bibitem[{Keilson \protect\BIBand{} Servi(1988)}]{NP:Keilson88}
Keilson J, Servi L (1988) A distributional form of {L}ittle's law.
  \emph{Operations Research Letters} 7(5):223--227.

\bibitem[{Kella \protect\BIBand{} Yechiali(1985)}]{NP:Kella85}
Kella O, Yechiali U (1985) Waiting times in the non-preemptive priority
  {M/M/$c$} queue. \emph{Communications in Statistics. Stochastic Models}
  1(2):257--262.

\bibitem[{Li \protect\BIBand{} Wong(2008)}]{NP:Li08}
Li Y, Wong R (2008) Integral and series representations of the {D}irac delta
  function. \emph{Communications in Pure and Applied Analysis} 7(2):229--247.

\bibitem[{Orman(1995)}]{NP:Orman95}
Orman A (1995) \emph{Models for Scheduling a Multifunction Phased Array Radar
  System}. Ph.D. thesis, University of Southampton, Southampton, UK.

\bibitem[{Orman et~al.(1996)Orman, Potts, Shahani, \protect\BIBand{}
  Moore}]{NP:Orman96}
Orman A, Potts C, Shahani A, Moore A (1996) Scheduling for a multifunction
  phased array radar system. \emph{European Journal of Operational Research}
  90(1):13--25.

\bibitem[{Pestalozzi(1964)}]{NP:Pestalozzi64}
Pestalozzi G (1964) Priority rules for runway use. \emph{Operations Research}
  12(6):941--950.

\bibitem[{Picard(1881)}]{NP:Picard81}
Picard E (1881) Sur une extension aux fonctions de deux variables du
  probl\`{e}me de {R}iemann relatif aux fonctions hyperg\'{e}om\'{e}triques.
  \emph{Annales Scientifiques de l'\'{E}cole Normale Sup\'{e}rieure, Serie 2}
  10:305--322.

\bibitem[{Press et~al.(2007)Press, Teukolsky, Vetterling, \protect\BIBand{}
  Flannery}]{NP:Press07}
Press W, Teukolsky S, Vetterling W, Flannery B (2007) \emph{Numerical Recipes:
  The Art of Scientific Computing} (Cambridge, UK: Cambridge University Press),
  third edition.

\bibitem[{Rice(1980)}]{NP:Rice80}
Rice S (1980) Distribution of quadratic forms in normal random variables --
  evaluation by numerical integration. \emph{SIAM Journal on Scientific and
  Statistical Computing} 1(4):438--448.

\bibitem[{Sleptchenko(2003)}]{NP:Sleptchenko03}
Sleptchenko A (2003) Multi-class, multi-server, queues with non-preemptive
  priorities. Report Eurandom Vol. 2003016, {EURANDOM}, Eindhoven University of
  Technology, The Netherlands,
  \urlprefix\url{https://pure.tue.nl/ws/portalfiles/portal/3290904/567759.pdf}.

\bibitem[{Taylor \protect\BIBand{} Templeton(1980)}]{NP:Taylor80}
Taylor I, Templeton J (1980) Waiting time in a multi-server cutoff-priority
  queue, and its application to an urban ambulance service. \emph{Operations
  Research} 28(5):1168--1188.

\bibitem[{Trefethen \protect\BIBand{} Weideman(2014)}]{NP:Trefethen14}
Trefethen L, Weideman J (2014) The exponentially convergent trapezoidal rule.
  \emph{SIAM Review} 56(3):385--458.

\bibitem[{Zuk \protect\BIBand{} Kirszenblat(2023{\natexlab{a}})}]{NP:Zuk23B}
Zuk J, Kirszenblat D (2023{\natexlab{a}}) Exact results for the distribution of
  the partial busy period for a multi-server queue,
  \urlprefix\url{http://dx.doi.org/10.2139/ssrn.4569933}, submitted for
  publication to \emph{Queueing Sytems}.

\bibitem[{Zuk \protect\BIBand{} Kirszenblat(2023{\natexlab{b}})}]{NP:Zuk23}
Zuk J, Kirszenblat D (2023{\natexlab{b}}) Explicit results for the
  distributions of queue lengths for a non-preemptive two-level priority queue,
  \urlprefix\url{http://dx.doi.org/10.2139/ssrn.4574550}, submitted for
  publication to \emph{Annals of Operations Research}.

\bibitem[{Zuk \protect\BIBand{} Kirszenblat(2023{\natexlab{c}})}]{NP:Zuk23A}
Zuk J, Kirszenblat D (2023{\natexlab{c}}) Joint queue-length distribution for
  the non-preemptive multi-server multi-level {M}arkovian priority queue,
  \urlprefix\url{http://dx.doi.org/10.2139/ssrn.4621771}, submitted for
  publication to \emph{Annals of Operations Research}.

\end{thebibliography}



\end{document}